# Some Applications of Arutyunov Mordukhovich Zhukovskiy Theorem to Stochastic Integral Equations

Jinlu Li

**Abstract**
Mordukhovich derivatives (Mordukhovich coderivatives) of set-valued mappings in Banach spaces have firmly laid the foundation of the theory of generalized differentiation in set-valued analysis, which has been widely applied to optimization theory, equilibrium theory, variational analysis, and so forth, with respect to set-valued mappings. One of the most important applications of Mordukhovich derivatives is to define the covering constants for set-valued mappings in Banach spaces, which play an important role in the well-known Arutyunov Mordukhovich Zhukovskiy Parameterized Coincidence Point Theorem (Theorem 3.1 in [1]). In [15], this theorem is simply named as AMZ Theorem. In this paper, we consider locally or globally stochastic infinitely dimensional systems of linear equations in $l_p$ space. We use the Mordukhovich derivatives to precisely find the covering constants for linear and continuous mappings in $l_p$ spaces. Then, by using the AMZ Theorem, we prove an existence theorem for solutions to some locally or globally stochastic infinitely dimensional systems of linear functional equations in $l_p$ spaces and an existence theorem for solutions to some stochastic integral equations.

1. Introduction

Let $(X, \|\cdot\|_X)$ and $(Y, \|\cdot\|_Y)$ be real Banach spaces with topological dual spaces $X^*$ and $Y^*$, and with origins $\theta_X$ and $\theta_Y$, respectively. Let $(S, \tau)$ be a topological space. Let $F(\cdot): X \rightrightarrows Y$ and $G(\cdot, \cdot): X \times S \rightrightarrows Y$ be set-valued mappings. In set-valued and variational analysis, the following parameterized coincidence point problems have attracted many authors' attention: Find an open subset $W \subset S$ and a single-valued mapping $\sigma: W \to X$ such that

$$F(\sigma(s)) \cap G(\sigma(s), s) \neq \emptyset, \text{ for any } s \in W. \tag{1.1}$$

An existence of solutions for the above parameterized coincidence point problems is proved in Theorem 3.1 in [1], which is called Arutyunov Mordukhovich Zhukovskiy Parameterized Coincidence Point Theorem. Considering the importance of this theorem, in [15] and in this paper, the Arutyunov Mordukhovich Zhukovskiy Parameterized Coincidence Point Theorem is simply named as AMZ Theorem. This theorem will be reviewed in section 2. The results of this theorem have been widely applied to set-valued analysis, such as set-valued optimization, set-valued equilibrium, set-valued variational inequality, and so forth (See [1–8, 16−22]). When we consider some specific different cases in (1.1), they can be applied to solving some corresponding problems, which are demonstrated by the following examples (1.3) to (1.5).

Let $A$ and $B$ be nonempty subsets of $Y$. We define $\|A - B\|_{\min} = \min\{\|v - w\|_Y : v \in A \text{ and } w \in B\}$. So, if $\|A - B\|_{\min}\ r$, for some $r \geq 0$, then, there are $a \in A$ and $b \in B$ such that

Department of Mathematics
Shawnee State University
Portsmouth, Ohio 45662 USA
Email: jli@shawnee.edu

**Keywords**  Fréchet differentiability; Mordukhovich derivative; covering constant; AMZ Theorem; locally and globally stochastic infinitely dimensional system of linear equations; stochastic integral equation

$$\|a - b\|_Y = \min\{\|v - w\|_Y : v \in A \text{ and } w \in B\} = r.$$

Then, the above parameterized coincidence points problem (1.1) becomes a parameterized minimization problem: Find an open subset $W \subset S$ and a single-valued mapping $\sigma: W \to X$ such that

$$\|F(\sigma(s)) - G(\sigma(s), s)\|_{\min} = 0, \text{ for any } s \in W. \tag{1.2}$$

This immediately implies that (1.1) is equivalent to the following parameterized equation problem: Find an open subset $W \subset S$ and a single-valued mapping $\sigma: W \to X$ such that

$$\theta_Y \in F(\sigma(s)) - G(\sigma(s), s), \text{ for any } s \in W. \tag{1.3}$$

In particular, if we study single-valued mappings, which are considered as special cases of set-valued mappings with values of singletons, then the parameterized coincidence point problems (1.1) has some special applications to parameterized (or stochastic) fixed point problems, parameterized (or stochastic) equations, and so forth.

Let $F(\cdot): X \to Y$ be a single-valued mapping and let $G(\cdot, \cdot): X \times S \rightrightarrows Y$ be set-valued mapping. The parameterized inclusion point problem with respect to $F$ and $G$ is to find an open subset $W \subset S$ and a single-valued mapping $\sigma: W \to X$ such that

$$F(\sigma(s)) \in G(\sigma(s), s), \text{ for any } s \in W. \tag{1.4}$$

For more special cases, let both $F: X \to Y$ and $G: X \times S \to Y$ be single-valued mappings. The parameterized equation problem with respect to $F$ and $G$ is to find an open subset $W \subset S$ and a single-valued mapping $\sigma: W \to X$ such that

$$F(\sigma(s)) = G(\sigma(s), s), \text{ for any } s \in W. \tag{1.5}$$

(1.4) and (1.5) have been used in [15] for studying the existence of solutions for some stochastic fixed-point problems. In this paper, we have more applications of the AMZ Theorem (1.1) and of the special cases (1.4) and (1.5) to stochastic integral equations.

This paper is organized as follows. In section 2, for the easy reference, we overview the Fréchet and Mordukhovich differentiability of mappings in Banach spaces and overview some concepts of covering constants for mappings in Banach spaces and the AMZ Theorem. We will give some special version for single valued mapping in the AMZ Theorem. In section 3, we study the solvability of some stochastic infinitely dimensional systems of linear functional equations in $l_p$ space. In section 4, we use the AMZ theorem to prove the existence of solutions to some stochastic integral equations in $L_2$ space.

## 2. Preliminaries

### 2.1. An Overview for Fréchet and Mordukhovich Differentiability in Banach Spaces

In this section, for easy reference, we first review the concepts of Gâteaux directional differentiability, Fréchet differentiability, and Strict Fréchet differentiability for single-valued mappings in Banach spaces. Then, we quickly review the concepts of Mordukhovich differentiability for set-valued mappings in Banach spaces. See [9–19, 22] for more details and [15,17] for a concise review.

Let $(X, \|\cdot\|_X)$ and $(Y, \|\cdot\|_Y)$ be real Banach spaces with topological dual spaces $X^*$ and $Y^*$, respectively. Let $\langle \cdot, \cdot \rangle_X$ denote the real canonical pairing between $X^*$ and $X$ and $\langle \cdot, \cdot \rangle_Y$ the real canonical pairing between $Y^*$ and $Y$. Let $\theta_X$ and $\theta_Y$ denote the origins in $X$ and $Y$, respectively. Let $\varphi: X \to Y$ be a single-valued mapping and let $\bar{x} \in X$. If there is a linear and continuous mapping $\nabla\varphi(\bar{x}): X \to Y$ such that

$$\lim_{x \to \bar{x}} \frac{\varphi(x) - \varphi(\bar{x}) - \nabla\varphi(\bar{x})(x - \bar{x})}{\|x - \bar{x}\|_X} = \theta_Y,$$

then $\varphi$ is said to be Fréchet differentiable at $\bar{x}$; $\nabla\varphi(\bar{x})$ is called the Fréchet derivative of $\varphi$ at $\bar{x}$.

The Mordukhovich derivative for set-valued mappings in Banach spaces forms the foundation of generalized differentiation in set-valued and variational analysis in Banach spaces (See [1–5, 16 – 19] Next, we review the concepts of Mordukhovich derivatives for set-valued mappings. See [16 –19] for more details. Let $\Delta$ be a nonempty subset of $X$ and let $F: \Delta \rightrightarrows Y$ be a set-valued mapping. The graph of $F$ is defined by the following subset in $\Delta \times Y$, $\mathrm{gph}F = \{(x, y) \in \Delta \times Y: y \in F(x)\}$. For $(x, y) \in \mathrm{gph}F$, that is, for $x \in \Delta$ and $y \in F(x)$, the Mordukhovich derivative of $F$ at point $(x, y)$ is a set valued mapping $\widehat{D}^*F(x, y): Y^* \rightrightarrows X^*$. For any $y^* \in Y^*$, it is defined by (see Definitions 1.13 and 1.32 in Chapter 1 in [18])

$$\widehat{D}^*F(x,y)(y^*) = \left\{ z^* \in X^*: \limsup_{\substack{(u,v) \to (x,y) \\ (u,v) \in \mathrm{gph}F}} \frac{\langle z^*, u-x \rangle_X - \langle y^*, v-y \rangle_Y}{\|u-x\|_X + \|v-y\|_Y} \leq 0 \right\}. \tag{2.1}$$

If $(x, y) \notin \mathrm{gph}F$, then, we define $\widehat{D}^*F(x, y)(y^*) = \emptyset$, for any $y^* \in Y^*$. The Mordukhovich derivative is also called Mordukhovich coderivative, Fréchet coderivative, or the coderivative of $F$. By the above definition (1.1), $\widehat{D}^*F(x, y): Y^* \rightrightarrows X^*$ is a set valued mapping, which also can be considered as the Mordukhovich differential operator of $F$ at point $(x, y)$. In particular, let $\varphi: \Delta \to Y$ be a single valued mapping. By (2.1), the Mordukhovich derivative of $\varphi$ at point $(x, \varphi(x))$ is a set-valued mapping $\widehat{D}^*\varphi(x, \varphi(x)): Y^* \rightrightarrows X^*$. For any $y^* \in Y^*$, denote $\widehat{D}^*\varphi(x, \varphi(x))(y^*)$ by $\widehat{D}^*\varphi(x)(y^*)$, which is defined by

$$\widehat{D}^*\varphi(x, \varphi(x))(y^*) = \widehat{D}^*\varphi(x)(y^*) = \left\{ z^* \in X^*: \limsup_{\substack{(u,\varphi(u)) \to (x,\varphi(x)) \\ u \in \Delta}} \frac{\langle z^*, u-x \rangle_X - \langle y^*, \varphi(u)-\varphi(x) \rangle_Y}{\|u-x\|_X + \|\varphi(u)-\varphi(x)\|_Y} \leq 0 \right\}.$$

Furthermore, if $\varphi: \Delta \to Y$ is a continuous single-valued mapping, then, for any $y^* \in Y^*$,

$$\widehat{D}^*\varphi(x, \varphi(x))(y^*) = \widehat{D}^*\varphi(x)(y^*) = \left\{ z^* \in X^*: \limsup_{\substack{u \to x \\ u \in \Delta}} \frac{\langle z^*, u-x \rangle_X - \langle y^*, \varphi(u)-\varphi(x) \rangle_Y}{\|u-x\|_X + \|\varphi(u)-\varphi(x)\|_Y} \leq 0 \right\}. \tag{2.2}$$

The following theorem shows the connection between Fréchet derivatives and Mordukhovich derivatives for single-valued mappings. Its results provide a powerful tool to calculate the Mordukhovich derivatives by the Fréchet derivatives of single-valued mappings.

**Theorem 1.38 in [18]**. *Let $X$ be a Banach space with dual space $X^*$ and let $\varphi: X \to Y$ be a single-valued mapping. Suppose that $f$ is Fréchet differentiable at $x \in X$ with $y = \varphi(x)$. Then, the Mordukhovich derivative of $\varphi$ at $x$ satisfies the following equation*

$$\widehat{D}^*\varphi(x, y)(y^*) = \{(\nabla\varphi(x))^*(y^*)\}, \text{ for all } y^* \in Y^*.$$

**Corollary 1.39 in [18].** (Coderivatives of linear and continuous operators $A$). *Let the conditions in Theorem* 1.38 *in* [18] *be fulfilled. Then*,

$$\widehat{D}^*A(x)(y^*) = \{A^*(y^*)\}, \text{ for all } y^* \in Y^*.$$

Mordukhovich derivatives have been widely applied to nonlinear analysis, such as operator theory, optimization theory, approximation theory, control theory, and so forth (see [1−5, 16 − 19]).

### 2.2. An Overview for Arutyunov Mordukhovich Zhukovskiy Theorem

One of the most important applications of the Mordukhovich derivatives of set-valued mappings is to define the covering constants for set-valued mappings. The existence of the covering constants for considered set-valued mappings is a sufficient condition in the AMZ Theorem (Theorem 3.1 in [1]). This theorem and some related results have played very important roles in set-valued analysis that has been widely applied to optimization theory, control theory, and so forth (See [1−5, 15−19]). To state the AMZ Theorem, we need to recall the concept of Asplund Banach spaces, the concepts of covering constants, and Lipschitz-like modulus for set-valued mappings in Banach spaces. For more details, see [1−5, 7, 20]. A Banach space $Z$ is Asplund if every convex continuous function defined on an open convex set $O$ in $Z$ is Fréchet differentiable on a dense subset of $O$. The class of Asplund Banach spaces is very large, which includes all reflexive Banach spaces (see [20], or page 182 in [1]). This implies that every uniformly convex and uniformly smooth Banach space is Asplund.

Let $X$ and $Y$ be Banach spaces and let $U$ and $V$ be nonempty subsets in $X$ and $Y$, respectively. Let $\Phi: X \rightrightarrows Y$ be a multifunction (a set-valued mapping). The graph of $\Phi$ is denoted by gph $\Phi$, which is a subset of $X \times Y$ defined by gph $\Phi \coloneqq \{(x,y) \in X \times Y: y \in \Phi(x)\}$. We say that $\Phi$ enjoys the covering property with modulus $\alpha > 0$ (or it has the $\alpha$-covering property) on $U$ relative to $V$ if (see (2.1) in [1])

$$\Phi(x) \cap V + \alpha r \mathbb{B}_Y \subset \Phi(x + r\mathbb{B}_X), \text{ whenever } x + r\mathbb{B}_X \subset U, \text{ as } r > 0. \qquad (2.3)$$

Here, $\mathbb{B}_X$ and $\mathbb{B}_Y$ are the closed unit balls in $X$ and $Y$, respectively. The supremum of all such moduli $\{\gamma\}$ in (2.3) is called the exact covering bound of $\Phi$ around $(x, y)$, which is denoted by

$$\text{cov}\Phi(x, y) = \sup\{\alpha: \alpha \text{ satisfies (2.3) for some } U \subset X, V \subset Y\}.$$

The multifunction $\Phi: X \rightrightarrows Y$ between $X$ and $Y$ is called Lipschitz-like on $U$ relative to $V$ with some modulus $\beta > 0$ if we have (see (2.3) in [1])

$$\Phi(x) \cap V \subset \Phi(u) + \beta \|x - u\|_X \mathbb{B}_Y, \text{ for all } x, u \in U. \qquad (2.4)$$

This is a natural extension of Lipschitz condition from single-valued mappings to set-valued mappings. In particular, let $\varphi: U \to Y$ be a single-valued mapping, $\varphi$ is said to satisfy the Lipschitz condition on $U$ relative to $V$ with respect to some modulus $\beta > 0$ if

$$\|\varphi(x) - \varphi(u)\|_Y \leq \beta \|x - u\|_X, \text{ for all } x, u \in U \text{ with } \varphi(x) \in V.$$

Note that, in [1], the covering property of set-valued mappings is defined in normed spaces. The covering property for set-valued mappings is a very important concept, which has been widely applied to game theory, optimization theory, equilibrium theory and variational analysis, with respect to set-valued mapping. For example, see [1−5, 16−19] for more details. One of the applications of Mordukhovich derivatives of set-valued mappings is to define the covering constants for set-valued mappings. The covering constant for $\Phi: X \rightrightarrows Y$ at point $(\bar{x}, \bar{y}) \in \text{gph } \Phi$ is defined by (see (2.6) in [1])

$$\hat{a}(\Phi, \bar{x}, \bar{y}) \coloneqq \sup_{\eta > 0} \inf\{\|z^*\|_{X^*}: z^* \in \widehat{D}^*\Phi(x, y)(w^*), x \in \mathbb{B}_X(\bar{x}, \eta), y \in \Phi(x) \cap \mathbb{B}_Y(\bar{y}, \eta), \|w^*\|_{Y^*} = 1\}. \qquad (2.5)$$

Here, $\|\cdot\|_{X^*}$ and $\|\cdot\|_{Y^*}$ denote the norms in $X^*$ and $Y^*$, respectively. $\mathbb{B}_X(\bar{x}, \eta)$ is the closed ball in $X$ centered at $\bar{x}$ with radius $\eta$, and $\mathbb{B}_Y(\bar{y}, \eta)$ is the closed ball in $Y$ centered at $\bar{y}$ with radius $\eta$.

In particular, let $\varphi: X \to Y$ be a single-valued mapping. For any $\bar{x}, \bar{y} \in X$ with $\bar{y} = \varphi(\bar{x})$, (2.5) becomes

$$\hat{\alpha}(\varphi, \bar{x}, \bar{y}) = \sup_{\eta > 0} \inf\{\|z^*\|_{X^*} : z^* \in \widehat{D}^*\varphi(x, y)(w^*), x \in \mathbb{B}_X(\bar{x}, \eta), y \in \mathbb{B}_Y(\bar{y}, \eta), \|w^*\|_{Y^*} = 1\}. \quad (2.6)$$

The following theorem precisely provides some more details about the connections between the local covering property and the covering constants for mappings around a given point.

**Theorem 4.1 in [18]** (neighborhood characterization of local covering) *Let $F: X \rightrightarrows Y$ be a set-valued mapping between Asplund spaces. Assume that F is closed-graph around $(\bar{x}, \bar{y}) \in$ gph $F$. Then the following are equivalent*:

(a) *F enjoys the local covering property around $(\bar{x}, \bar{y})$ (that is, cov $F(\bar{x}, \bar{y}) > 0$).*
(b) *One has $\hat{\alpha}(F, \bar{x}, \bar{y}) > 0$.*

*Moreover, the exact covering bound of F around $(\bar{x}, \bar{y})$ is computed by*

$$\text{cov} F(\bar{x}, \bar{y}) = \hat{\alpha}(F, \bar{x}, \bar{y}).$$

In particular, the exact covering bounds of Fréchet differentiable single-valued mappings have been studied in [18]. For example, we have the following results from [18], which can be used in Sections 3, 4.

**Corollary 1.58 in [18]** (covering for linear operators) *A linear and continuous operator $A: X \to Y$ has the covering property around point $\bar{x} \in X$ with $\bar{y} = A(\bar{x})$ if and only if A is surjective and we have*

$$\text{cov} A(\bar{x}, \bar{y}) = \inf\{\|A^* y^*\|_{X^*} : \|y^*\|_{Y^*} = 1\}, \text{ for all } \bar{x} \in X.$$

Let $F: X \rightrightarrows Y$ be a set valued mapping. Let $\bar{x} \in X$ and $\bar{y} \in Y$ with $\bar{y} \in F(\bar{x})$. The multifunction $F$ is said to be closed around $(\bar{x}, \bar{y})$ if there exist neighborhood $U$ of $\bar{x}$ and $V$ of $\bar{y}$ such that (gph $F$)$\cap$ (cl$U \times$ cl$V$) is closed in $X \times Y$. Where, "cl" indicates the topological closure operation. In particular, let $F: X \to Y$ be a single-valued mapping. Let $\bar{x} \in X$ and $\bar{y} \in Y$ with $\bar{y} = F(\bar{x})$. Then, the continuity of $F$ at $(\bar{x}, \bar{y})$ implies that $F$ is closed around $(\bar{x}, \bar{y})$. Next, we state the AMZ Theorem (Theorem 3.1 in [1])

**(AMZ Theorem)** *Let the Banach spaces $X$ and $Y$ in be Asplund and let $P$ be a topological space. Let $F: X \rightrightarrows Y$ and $G(\cdot, \cdot): X \times P \rightrightarrows Y$ be set-valued mappings. Let $\bar{x} \in X$ and $\bar{y} \in Y$ with $\bar{y} \in F(\bar{x})$. Suppose that the following conditions are satisfied:*

(A1) *The multifunction $F: X \rightrightarrows Y$ is closed around $(\bar{x}, \bar{y})$.*

(A2) *There are neighborhoods $U \subset X$ of $\bar{x}$, $V \subset Y$ of $\bar{y}$, and $O$ of $\bar{p} \in P$ as well as a number $\beta \geq 0$ such that the multifunction $G(\cdot, p): X \rightrightarrows Y$ is Lipschitz-like on $U$ relative to $V$ for each $p \in O$ with the uniform modulus $\beta$, while the multifunction $p \to G(\bar{x}, p)$ is lower/inner semicontinuous at $\bar{p}$.*

(A3) *The Lipschitzian modulus $\beta$ of $G(\cdot, p)$ is chosen as $\beta < \hat{\alpha}(F, \bar{x}, \bar{y})$, where $\hat{\alpha}(F, \bar{x}, \bar{y})$ is the covering constant of F around $(\bar{x}, \bar{y})$ taken from (2.5).*

*Then for each $\alpha > 0$ with $\beta < \alpha < \hat{\alpha}(F, \bar{x}, \bar{y})$, there exist a neighborhood $W \subset P$ of $\bar{p}$ and a single-valued mapping $\sigma: W \to X$ such that whenever $p \in W$ we have*

$$F(\sigma(p)) \cap G(\sigma(p), p) \neq \emptyset \quad \text{and} \quad \|\sigma(p) - \bar{x}\|_X \leq \frac{\text{dist}(\bar{y}, G(\bar{x}, p))}{\alpha - \beta}. \quad (2.7)$$

This result of Theorem 3.1 in [1] is very powerful. We provide some immediate consequences of the

AMZ Theorem for some special cases, in which at least one of $F$ and $G$ in AMZ Theorem is a single-valued mapping. In the first corollary below, we let the mapping $F$ be single-valued and $G$ be set-valued.

**Corollary 2.1.** *Let the Banach spaces $X$ and $Y$ in be Asplund and let $P$ be a topological space. Let $F: X \to Y$ be a single-valued mapping and let $G(\cdot, \cdot): X \times P \rightrightarrows Y$ be a set-valued mapping. Let $\bar{x} \in X$ and $\bar{y} \in Y$ with $\bar{y} = F(\bar{x})$. Suppose that the following conditions are satisfied:*

(A1) *The mapping $F: X \to Y$ is continuous around $(\bar{x}, \bar{y})$.*

(A2) *There are neighborhoods $U \subset X$ of $\bar{x}$, $V \subset Y$ of $\bar{y}$, and $O$ of $\bar{p} \in P$ as well as a number $\beta \geq 0$ such that the multifunction $G(\cdot, p): X \rightrightarrows Y$ is Lipschitz-like on $U$ relative to $V$ for each $p \in O$ with the uniform modulus $\beta$, while the multifunction $p \to G(\bar{x}, p)$ is lower/inner semicontinuous at $\bar{p}$.*

(A3) *The Lipschitzian modulus $\beta$ of $G(\cdot, p)$ is chosen as $\beta < \hat{\alpha}(F, \bar{x}, \bar{y})$, where $\hat{\alpha}(F, \bar{x}, \bar{y})$ is the covering constant of $F$ around $(\bar{x}, \bar{y})$ taken from (2.5).*

*Then for each $\alpha > 0$ with $\beta < \alpha < \hat{\alpha}(F, \bar{x}, \bar{y})$, there exist a neighborhood $W \subset P$ of $\bar{p}$ and a single-valued mapping $\sigma: W \to X$ such that whenever $p \in W$ we have*

$$F(\sigma(p)) \in G(\sigma(p), p) \quad \text{and} \quad \|\sigma(p) - \bar{x}\|_X \leq \frac{\text{dist}(\bar{y}, G(\bar{x}, p))}{\alpha - \beta}. \tag{2.8}$$

*Proof.* The proof of this corollary is straight forward and it is omitted here. □

In next corollary, we let both the considered mapping $F$ and $G$ in the AMZ Theorem be single-valued.

**Corollary 2.2.** *Let the Banach spaces $X$ and $Y$ be Asplund and let $P$ be a topological space. Let $F: X \to Y$ and $G(\cdot, \cdot): X \times P \to Y$ be single-valued mappings. Let $\bar{x} \in X$ and $\bar{y} \in Y$ with $\bar{y} = F(\bar{x})$. Suppose that the following conditions are satisfied:*

(A1) *The mapping $F: X \to Y$ is continuous around $(\bar{x}, \bar{y})$.*

(A2) *There are neighborhoods $U \subset X$ of $\bar{x}$, $V \subset Y$ of $\bar{y}$, and $O$ of $\bar{p} \in P$ as well as a number $\beta \geq 0$ such that the mapping $G(\cdot, p): X \to Y$ satisfies the Lipschitz condition on $U$ relative to $V$ for each $p \in O$ with the uniform modulus $\beta$, while the mapping $p \to G(\bar{x}, p)$ is lower semicontinuous at $\bar{p}$.*

(A3) *The Lipschitzian modulus $\beta$ of $G(\cdot, p)$ is chosen as $\beta < \hat{\alpha}(F, \bar{x}, \bar{y})$, where $\hat{\alpha}(F, \bar{x}, \bar{y})$ is the covering constant of $F$ around $(\bar{x}, \bar{y})$ taken from (2.5).*

*Then for each $\alpha > 0$ with $\beta < \alpha < \hat{\alpha}(F, \bar{x}, \bar{y})$, there exist a neighborhood $W \subset P$ of $\bar{p}$ and a single-valued mapping $\sigma: W \to X$ such that whenever $p \in W$ we have*

$$F(\sigma(p)) = G(\sigma(p), p) \quad \text{and} \quad \|\sigma(p) - \bar{x}\|_X \leq \frac{\|G(\bar{x}, p) - \bar{y}\|_Y}{\alpha - \beta}. \tag{2.9}$$

*Proof.* The proof of this corollary is straight forward and it is omitted here. □

**Corollary 2.3.** *Let the Banach spaces $X$ and $Y$ in be Asplund and let $P$ be a topological space. Let $F: X \rightrightarrows Y$ be a set-valued mapping and let $G(\cdot, \cdot): X \times P \to Y$ be a set-valued mapping. Let $\bar{x} \in X$ and $\bar{y} \in Y$ with $\bar{y} = F(\bar{x})$. Suppose that the following conditions are satisfied:*

(A1) *The multifunction $F: X \rightrightarrows Y$ is closed around $(\bar{x}, \bar{y})$.*

(A2) *There are neighborhoods $U \subset X$ of $\bar{x}$, $V \subset Y$ of $\bar{y}$, and $O$ of $\bar{p} \in P$ as well as a number $\beta \geq 0$ such that the mapping $G(\cdot, p): X \to Y$ satisfies the Lipschitz condition on $U$ relative to $V$ for each $p \in O$ with the uniform modulus $\beta$, while the mapping $p \to G(\bar{x}, p)$ is lower semicontinuous at $\bar{p}$.*

(A3) *The Lipschitzian modulus $\beta$ of $G(\cdot, p)$ is chosen as $\beta < \hat{\alpha}(F, \bar{x}, \bar{y})$, where $\hat{\alpha}(F, \bar{x}, \bar{y})$ is the covering constant of $F$ around $(\bar{x}, \bar{y})$ taken from (2.5).*

*Then for each $\alpha > 0$ with $\beta < \alpha < \hat{\alpha}(F, \bar{x}, \bar{y})$, there exist a neighborhood $W \subset P$ of $\bar{p}$ and a single-valued mapping $\sigma: W \to X$ such that whenever $p \in W$ we have*

$$G(\sigma(p), p) \in F(\sigma(p)) \quad \text{and} \quad \|\sigma(p) - \bar{x}\|_X \leq \frac{\text{dist}(\bar{y}, G(\bar{x}, p))}{\alpha - \beta}.$$

*Proof.* The proof of this corollary is straight forward and it is omitted here. □

The results of the AMZ Theorem are very strong, which provides a general and powerful tool to prove some existence problems in nonlinear analysis. The themes of this paper are to prove some stochastic integral equations by the AMZ Theorem.

However, we note that, for a given set-valued mapping $F$, it is very difficult to calculate Mordukhovich derivatives of $F$ (see (2.1)). Even for single-valued mapping $F$, it is still very difficult to calculate the Mordukhovich derivatives of $F$ (see (2.2)), except some special cases (see Corollary 1.39 in [18]). One more step further, since the covering constants for both set-valued and single-valued mappings are defined by its Mordukhovich derivatives, one knows immediately that it is extremally difficult to calculate the covering constants for the considered set-valued or single-valued mappings in the AMZ Theorem, except some special cases (see Corollary 1.58 in [18]).

By Theorem 4.1 in [18], for the considered mapping $F$ under some conditions, the covering constants $\hat{\alpha}(F, \bar{x}, \bar{y})$ for $F$ at a point $(\bar{x}, \bar{y})$ can be calculated by the exact covering bound $\text{cov}F(\bar{x}, \bar{y})$ of $F$ at $(\bar{x}, \bar{y})$. Meanwhile, the calculation for $\text{cov}F(\bar{x}, \bar{y})$ is also very complicated, in general (see (2.3)., which can be demonstrated by the following example.

**Example 2.1** (Example 2 in [4] and Example 4.2 in [2]). Let $\theta$ denote the origin of Euclidean space $\mathbb{R}^2$. Define a single-valued mapping $F: \mathbb{R}^2 \to \mathbb{R}^2$ by

$$F((x_1, x_2)) = \left( \frac{x_1^2 - x_2^2}{\sqrt{x_1^2 + x_2^2}}, \frac{2x_1 x_2}{\sqrt{x_1^2 + x_2^2}} \right), \text{ for } (x_1, x_2) \in \mathbb{R}^2 \setminus \{\theta\} \text{ with } F(\theta) = \theta.$$

Then, $F$ is continuous on $\mathbb{R}^2$. In [2, 4], the authors gave an elegant proof for the following result:

$$\alpha(F, \theta, \theta) = \hat{\alpha}(F, \theta, \theta) = 1.$$

The proof of the first equation $\alpha(F, \theta, \theta) = \hat{\alpha}(F, \theta, \theta)$ is based on Theorem 4.1 in [18], in which the Mordukhovich derivative of $F$ at $(\theta, \theta)$ is not used. One sees that in this example, the underlying space is the Euclidean space $\mathbb{R}^2$ and the considered single-valued mapping $F: \mathbb{R}^2 \to \mathbb{R}^2$ is not complicated. But the proof of the results that $\alpha(F, \theta, \theta) = 1$ in [4] is elegant and not easy.

Sometimes, the results of the covering constants for some mappings are very peculiar. It is well-known that the standard metric projection operator is an extremally important mapping in approximation theory, fixed point theory, optimization theory, and so forth. The following results of the covering constants for the standard metric projection operator are very surprising.

**Theorem 3.1 in [13].** *Let $(X, \|\cdot\|)$ be a real uniformly convex and uniformly smooth Banach space and let $\mathbb{B}$ denote the unit closed ball in $X$ with topological interior $\mathbb{B}^o$. For any $r > 0$, let $P_{r\mathbb{B}}: X \to r\mathbb{B}$ be the standard metric projection. For $\bar{x} \in X$ with $\bar{y} = P_{r\mathbb{B}}(\bar{x})$, the covering constant for the metric projection $P_{r\mathbb{B}}$ at $(\bar{x}, \bar{y})$ satisfies*

(a) $\hat{\alpha}(P_{r\mathbb{B}}, \bar{x}, \bar{y}) = 1$, *for any $\bar{x} \in r\mathbb{B}^o$*;
(b) $\hat{\alpha}(P_{r\mathbb{B}}, \bar{x}, \bar{y}) = 0$, *for any $\bar{x} \in X \backslash (r\mathbb{B}^o)$.*

By the singularity of the result of part (b) in the above theorem, we realize that the standard metric projection operator cannot be considered (as the mapping $F$) in the applications of the AMZ Theorem.

When we apply the AMZ Theorem to prove some existence problems, since the difficulty for finding the covering constants of the considered mappings (it may be impossible to find it), in general, this is why, in this paper, we only study some special mappings, for which the covering constants can be calculated.

## 3. Stochastic Systems of Linear Equations in $l_p$ Spaces

### 3.1 Operator Norms and Covering Constants of Linear and Continuous Mappings in $l_p$ Spaces

Let $p$, $q$ be positive numbers satisfying $1 < p, q < \infty$ and $\frac{1}{p} + \frac{1}{q} = 1$. $(l_p, \|\cdot\|_p)$ and $(l_q, \|\cdot\|_q)$ denote the standard real uniformly convex and uniformly smooth Banach spaces of sequences of real numbers, which are dual spaces to each other with the real canonical product $\langle \cdot, \cdot \rangle$ between $l_q$ and $l_p$. The origins of both $l_p$ and $l_q$ are exactly the same $\theta = \theta^* = (0, 0, \dots)$ Let $T$ denote the collection of all sequences of real numbers. Let $\mathbb{B}_p$ and $\mathbb{S}_p$ denote the closed unit ball and unit sphere in $l_p$, respectively. For any $x \in l_p$ and $r > 0$, let $\mathbb{B}_p(x, r)$ and $\mathbb{S}_p(x, r)$ respectively denote the closed ball and sphere in $l_p$ with center $x$ and radius $r$.

Let $A = (a_{ij})_{i,j=1}^{\infty}$ be a real $\infty \times \infty$ matrix (It is also named by a double sequence of real numbers). Under certain conditions on $A$, this real $\infty \times \infty$ matrix $A$ defines a linear mapping from $l_p$ to the set of all sequences of real numbers, such that

$$A(x) = (x_1, x_2, \dots) A = \left(\sum_{i=1}^{n} a_{ij} x_i\right)_{j=1}^{\infty}, \text{ for any } x = (x_1, x_2, \dots) \in l_p. \quad (3.1)$$

Where, $x = (x_1, x_2, \dots)$ and $A(x) = \left(\sum_{i=1}^{n} a_{ij} x_i\right)_{j=1}^{\infty} = \left(\sum_{i=1}^{n} a_{i1} x_i, \sum_{i=1}^{n} a_{i2} x_i, \dots\right)$ are also considered as $1 \times \infty$ matrices. In (3.1), for the mapping $A$ induced by matrix $A$, $A(x)$ is the value of the mapping $A$ at the point $x = (x_1, x_2, \dots) \in l_p$. The value is defined by $(x_1, x_2, \dots) A$ as the product of an $1 \times \infty$ matrix and an $\infty \times \infty$ matrix. It is well defined on $l_p$ if $A$ satisfies some conditions (see the following lemmas). We have some notations. Let $A = (a_{ij})_{i,j=1}^{\infty}$ be a real $\infty \times \infty$ matrix, if the corresponding mapping $A: l_p \to l_p$ is a linear and continuous mapping, then, we let $\|A\|_{\text{op}}$ denote the operator norm of $A$. We write

$$\|A\|_{\text{op}} = \sup \{\|A(x)\|_p : x \in \mathbb{S}_p\} \quad \text{and} \quad \|A\|_{\text{inf}} = \inf \{\|A(x)\|_p : x \in \mathbb{S}_p\}.$$

Notice that both $\|A\|_{\text{op}}$ and $\|A\|_{\text{inf}}$ depend on $p$.

**Lemma 3.1.** *Let $A = (a_{ij})_{i,j=1}^{\infty}$ be a real $\infty \times \infty$ matrix. Suppose that $A$ satisfies the following conditions*

$$\sum_{j=1}^{\infty}\left(\sum_{i=1}^{\infty}|a_{ij}|^q\right)^{\frac{p}{q}} < \infty \quad \text{and} \quad \sum_{i=1}^{\infty}\left(\sum_{j=1}^{\infty}|a_{ij}|^p\right)^{\frac{q}{p}} < \infty. \quad (3.2)$$

*Then, we have*

(i) *A defines a linear and continuous mapping $A: l_p \to l_p$ with operator norm $\|A\|_{op}$ satisfying*

$$\|A\|_{op} \leq \left(\sum_{j=1}^{\infty}\left(\sum_{i=1}^{\infty}|a_{ij}|^q\right)^{\frac{p}{q}}\right)^{\frac{1}{p}}; \tag{3.3}$$

(ii) *The adjoint operator $A^*$ of $A$ is its transpose $A^T$ that is also a linear and continuous mapping $A^T: l_q \to l_q$ with operator norm $\|A^T\|_{op}$ satisfying*

$$\|A^T\|_{op} \leq \left(\sum_{i=1}^{\infty}\left(\sum_{j=1}^{\infty}|a_{ij}|^p\right)^{\frac{q}{p}}\right)^{\frac{1}{q}}. \tag{3.4}$$

*Proof.* Proof of part (i). For any $x = (x_1, x_2, \ldots) \in l_p$, by (3.1), (3.2) and by Hölder inequality, we have

$$\|A(x)\|_p = \left(\sum_{j=1}^{\infty}\left|\sum_{i=1}^{\infty}a_{ij}x_i\right|^p\right)^{\frac{1}{p}} \leq \left(\sum_{j=1}^{\infty}\left(\sum_{i=1}^{\infty}|a_{ij}|^q\right)^{\frac{p}{q}}\right)^{\frac{1}{p}}\|x\|_p, \text{ for any } x = (x_1, x_2, \ldots) \in l_p.$$

This implies (3.3). Part (ii) can be similarly proved. Similarly, to (3.3) for each $q > 1$, we have

$$\|A^T(y)\|_q = \left(\sum_{i=1}^{\infty}\left|\sum_{j=1}^{\infty}a_{ij}y_j\right|^q\right)^{\frac{1}{q}} \leq \left(\sum_{i=1}^{\infty}\left(\sum_{j=1}^{\infty}|a_{ij}|^p\right)^{\frac{q}{p}}\right)^{\frac{1}{q}}\|y\|_q, \text{ For any } y = (y_1, y_2, \ldots) \in l_q.$$

This implies (3.4). □

**Theorem 3.2.** *Let $A = (a_{ij})_{i,j=1}^{\infty}$ be a real $\infty \times \infty$ matrix. Suppose that $A$ satisfies conditions (3.2). Then, we have*

(i) *$A$ is Fréchet differentiable at every point in $l_p$ with $\nabla(A)(x) = A$, for any $x \in l_p$;*

(ii) *The Mordukhovich derivative of $A$ satisfies that $\widehat{D}^*(A)(x, A(x)) = A^T$, for any $x \in l_p$;*

(iii) *Suppose that $\sum_{i=1}^{\infty}\left(\sum_{j=1}^{\infty}|a_{ij}|^p\right)^{\frac{q}{p}} \leq 1$, then the covering constant for $A$ is constant in $l_p$ with*

$$\hat{\alpha}(A, x, A(x)) = \|A^T\|_{\inf} = \inf\{\|A^T(y)\|_q: y \in \mathbb{S}_q\} = 0, \text{ for any } x \in l_p.$$

*Proof.* Proof of Part (i). By Lemma 3.1, the mapping $A: l_p \to l_p$ defined by (3.1) is a linear and continuous single-valued mapping. We have

$$\lim_{u \to x}\frac{A(u) - A(x) - A(u-x)}{\|u-x\|_p} = \theta, \text{ for any given } x \in l_p.$$

By Definition 1.13 in [17], this proves (i). By Theorem 1.38 in [17], part (i) induces part (ii) immediately. Now, by (ii), we prove part (iii). Since $A$ satisfies $\sum_{j=1}^{\infty}\left(\sum_{i=1}^{\infty}|a_{ij}|^q\right)^{\frac{p}{q}} \leq 1$, by (3.3), we have

$$\|A\|_{\mathrm{op}} \leq \left(\sum_{j=1}^{\infty}\left(\sum_{i=1}^{\infty}|a_{ij}|^q\right)^{\frac{p}{q}}\right)^{\frac{1}{p}} \leq 1. \tag{3.5}$$

Then, for any $x \in l_p$ and for any $\eta > 0$, similarly, by condition (3.5), we have

$$u \in \mathbb{B}_p(x,\eta) \implies A(u) \in \mathbb{B}_p(A(x),\eta), \text{ for any } u \in l_p. \tag{3.6}$$

For any $m \geq 1$, let $s_m$ denote the sequence of real numbers that has $m^{\text{th}}$ entry 1 and all other entries 0. By the condition (3.2), we have that

$$\sum_{j=1}^{\infty}\left(\sum_{i=1}^{\infty}|a_{ij}|^q\right)^{\frac{p}{q}} < \infty \implies \lim_{i \to \infty}\left(\sum_{i=1}^{\infty}|a_{ij}|^q\right)^{\frac{p}{q}} = 0.$$

Then, by (3.6) and (ii) and the above property of $A$, we calculate the covering constant for $A$ at an arbitrarily given point $x \in l_p$.

$\hat{\alpha}(A, x, A(x))$

$= \sup\limits_{\eta>0} \inf\{\|w\|_q : w \in \widehat{D}^*(A)(u, A(u))(y), u \in \mathbb{B}_p(x,\eta), A(u) \in \mathbb{B}_p(A(x),\eta), \|y\|_q = 1\}$

$= \sup\limits_{\eta>0} \inf\{\|A^T(y)\|_q : \{A^T(y)\} = \widehat{D}^*(A)(u, A(u))(y), u \in \mathbb{B}_p(x,\eta), A(u) \in \mathbb{B}_p(A(x),\eta), \|y\|_q = 1\}$

$= \sup\limits_{\eta>0} \inf\{\|A^T(y)\|_q : y \in l_q, \|y\|_q = 1\}$

$\leq \inf\{\|A^T(y)\|_q : y \in l_q, \|y\|_q = 1\}$

$\leq \sup\limits_{\eta>0} \inf\{\|A^T(s_m)\|_q : s_m \in l_q, \|s_m\|_q = 1, m = 1, 2, \ldots\}$

$\leq \inf\{\|A^T(s_m)\|_q : s_m \in l_q, \|s_m\|_q = 1, m = 1, 2, \ldots\}$

$= \inf\left\{(\sum_{i=1}^{\infty}|a_{im}|^q)^{\frac{1}{q}} : s_m \in l_q, \|s_m\|_q = 1, m = 1, 2, \ldots\right\}$

$= 0.$  □

We see that since $A^T : l_q \to l_q$ is a linear and continuous mapping, then, Part (iii) of Theorem 3.2 can be proved by using Corollary 1.58 in [18] (covering for linear and continuous operators).

Notice that part (iii) follows from the covering criterion (see Theorem 2.1 in [1]) and well-known facts of classical linear analysis (without any additional assumptions on $\|A\|$ or $\|A^T\|$), see, for instance, Lemma 1.18 in [18].

In particular, let $n$ be a positive integer and let $A = (a_{ij})_{i,j=1}^n$ be a real $n \times n$ matrix. In [15], Li proved that if $\det(A) \neq 0$ and $\sum_{i,j=1}^n a_{ij}^2 \leq 1$, then the covering constant for $A$ is constant in $\mathbb{R}^n$ satisfying

$$0 < \hat{\alpha}(A, x, A(x)) \leq \|A^T\|_{\mathrm{op}} \leq \left(\sum_{i,j=1}^n a_{ij}^2\right)^{\frac{1}{2}} \leq 1, \text{ for any } x \in \mathbb{R}^n. \tag{3.7}$$

This is proved based on that the unit sphere of $\mathbb{R}^n$ is a compact subset in $\mathbb{R}^n$. Let $\theta_n$ denote the origin of $\mathbb{R}^n$. We know that $\det(A) \neq 0$ if and only if the following system of linear equations has only solution $\theta_n$

$$xA = \theta_n. \tag{3.8}$$

It is equivalent to say that if the system of linear equations (3.8) has only solution $\theta_n$ and $\sum_{i,j=1}^{n} a_{ij}^2 \leq 1$, then the covering constant for $A$ is a positive constant in $\mathbb{R}^n$ satisfying (3.7). But, if $A = (a_{ij})_{i,j=1}^{\infty}$ is a real $\infty \times \infty$ matrix, then $\det(A)$ is undefined; and therefore, (3.8) cannot be defined by the undefined condition $\det(A) \neq 0$. Next, we provide a simple counterexample to show that, in $l_p$, the condition (3.8) that the system has only solution $\theta$ does not assure that the covering constant for $A$ is positive.

**Example 3.3.** Consider Hilbert space $l_2$. Let $A = (a_{ij})_{i,j=1}^{\infty}$ be a real $\infty \times \infty$ diagonal matrix with $a_{ii} = \frac{1}{i+1}$, for $i = 1, 2, \ldots$. Then, $A$ satisfies that, for any $x = (x_1, x_2, \ldots) \in l_2$, $xA = \theta$, if and only if $x = \theta$.

However, $\hat{\alpha}(A, x, A(x)) = 0$, for any $x = (x_1, x_2, x_3, \ldots) \in l_2$.

*Proof.* The mapping $A$ on $l_2$ induced by this diagonal matrix $A = (a_{ij})_{i,j=1}^{\infty}$ is a pointwise multiplication operator on $l_2$. That is, for any $x = (x_1, x_2, x_3, \ldots) \in l_2$, $A(x) = xA = (\frac{1}{2}x_1, \frac{1}{3}x_2, \frac{1}{4}x_3, \ldots) \in l_2$. It satisfies $A^* = A^T = A$. We can check that $A$ satisfies conditions (3.2) with respect to $p = q = 2$:

$$\sum_{j=1}^{\infty} \frac{1}{(j+1)^2} < 1 \quad \text{and} \quad \sum_{i=1}^{\infty} \frac{1}{(i+1)^2} < 1.$$

Now we show that $\hat{\alpha}(A, x, A(x)) = 0$, for any $x = (x_1, x_2, \ldots) \in l_2$. To this end, we calculate $\hat{\alpha}(A, x, A(x))$, with respect to an arbitrarily given $x = (x_1, x_2, \ldots) \in l_2$. For any $j = 1, 2, \ldots$, let $y^{(j)} \in l_2$, in which, the $j$th entry is 1 and all other entries are 0. We have

$$A^T(y^{(j)}) = y^{(j)}A^T = \frac{1}{j+1}y^{(j)}, \text{ for } j = 1, 2, \ldots. \tag{3.9}$$

This implies

$$\|A^T(y^{(j)})\|_2 = \left(\frac{1}{(j+1)^2}\right)^{\frac{1}{2}} = \frac{1}{j+1}, \text{ for } j = 1, 2, \ldots. \tag{3.10}$$

For any $x \in l_2$, by part (iii) of Theorem 3.2 and by (3.9) and (3.10), we have

$$\hat{\alpha}(A, x, A(x)) = \inf\{\|A^T(y)\|_2 : y \in l_2, \|y\|_2 = 1\}$$

$$\leq \inf\left\{\|A^T(y^{(j)})\|_2 : j = 1, 2, \ldots\right\} \leq \inf\left\{\frac{1}{j+1} : j = 1, 2, \ldots\right\} = 0. \qquad \square$$

### 3.2 Stochastic Systems of Linear Equations in $l_p$ Spaces

Let $(S, \tau, \mu)$ be a topological probability space in which $S$ is the sample space such that the topology $\tau$ on $S$ coincides with the $\sigma$-field of all events in $S$, and $\mu$ is the probability measure in $S$ defined on $\tau$. Let $A = (a_{ij})_{i,j=1}^{\infty}$ and $B = (b_{ij})_{i,j=1}^{\infty}$ be real $\infty \times \infty$ matrices. Suppose that $A$ and $B$ respectively define linear and continuous mappings from $l_p$ to itself. Let $\omega: S \to l_p$ be a $\tau$-measurable single-valued mapping (that is considered as a noise). Let $\bar{s} \in S$. If there is an event (neighborhood) $W \subset S$ of $\bar{s}$ and a single-valued mapping $\sigma: W \to l_p$ such that

$$A(\sigma(s)) = B(\sigma(s)) + \omega(s), \text{ for any } s \in W, \tag{3.11}$$

then $\sigma$ is called a solution to the *locally stochastic infinitely dimensional system of linear equations on W around the possible outcome $\bar{s}$ and with respect to the linear and continuous mappings A, B and the noise $\omega$*. In particular, in the inclusion property (2.4), if the neighborhood of $\bar{s}$ is the whole space $S$; that is, if

$$A(\sigma(s)) = B(\sigma(s)) + \omega(s), \text{ for any } s \in S,$$

then $\sigma$ is called a solution to the *globally stochastic infinitely dimensional system of linear equations on S around the possible outcome $\bar{s}$ and with respect to the linear and continuous mappings A, B and the noise $\omega$*. In order to study the solvability of the locally or globally stochastic systems of linear equations with respect to some given $\infty \times \infty$ matrices $A$ and $B$, we need to consider some conditions for these matrices and their dual matrices.

In this subsection, we use the AMZ Theorem to study the solubility of the locally or globally stochastic infinitely dimensional system of linear equations on a topological probability space $S$ around any given possible outcome in $S$ and with respect to some given matrices $A$, $B$ and a noise $\omega$. More precisely speaking, we will prove the existence of solutions to the locally stochastic systems of infinitely dimensional linear equations on a topological probability space $S$ around any given possible outcome $\bar{s}$ and with respect to $\infty \times \infty$ matrices $A$, $B$ and a noise $\omega$ with values in $l_p$.

Notice that the results of AMZ Theorem contain two parts:

(a) The first part is the solution $\sigma$ existence for a considered problem. But this theorem does not make sure whether the solution $\sigma$ is measurable;
(b) The second part provides an estimation of the difference between the solution $\sigma$ and the started point.

We believe that in the further study, the second part can be used to investigate the properties of the solution $\sigma$, which include the measurability. Once the measurability of $\sigma$ is proved, then $\sigma$ will automatically become a random variable defined on the underlying probability space $S$. This is why the problem (3.1) is preceded as a problem of stochastic infinitely dimensional system of linear equations and $\sigma$ is dubbed as a stochastic solution to problem (3.1).

One sees that the following Theorem may be proved by other way without using the AMZ Theorem. However, based on the themes of this paper, we will use the AMZ Theorem to prove the following theorem. As we mentioned in Section 2, the AMZ Theorem is a very important theorem in nonlinear analysis and it is not easy to be used to solve some specific problems. We consider the proof of the following theorem being attempting to use the AMZ Theorem to prove the solution existence for some specific problems.

We believe that with the development of generalized differentiation in set-valued analysis, some practical and feasible techniques for calculating the Mordukhovich derivatives and the covering constants for both set-valued and single-valued mappings will be obtained. Then, the AMZ Theorem will be used to prove the solution existence for more general problems with respect to more complicated mappings.

**Theorem 3.4.** *Let $(S, \tau, \mu)$ be a topological probability space. Let $A = (a_{ij})_{i,j=1}^{\infty}$ and $B = (b_{ij})_{i,j=1}^{\infty}$ be real $\infty \times \infty$ matrices. Suppose that A and B respectively define linear and continuous mappings from $l_p$ to itself. Let $\omega: S \to l_p$ be a $\tau$-measurable single-valued mapping. Let $\bar{s} \in S$. Suppose that $A, B$ and $\omega$ satisfy the following conditions.*

(a1) $0 < \|B\|_{op} < \hat{\alpha}(A, x, A(x)) \leq \|A^*\|_{op} \leq 1$, *for every* $x \in l_p$;

(a2) *The function* $s \to \omega(s)$ *is lower semicontinuous at* $\bar{s}$.

*Then, for any* $\lambda, \alpha$ *with* $\|B\|_{op} < \alpha < \lambda \leq \hat{\alpha}(A, x, A(x))$, *there exist a neighborhood* $W_{\lambda\alpha} \subset S$ *of* $\bar{s}$ *and a single-valued mapping* $\sigma_{\lambda\alpha}: W_{\lambda\alpha} \to l_p$ *such that*

$$A(\sigma_{\lambda\alpha}(s)) = B(\sigma_{\lambda\alpha}(s)) + \omega(s), \text{ for every } s \in W_{\lambda\alpha}, \tag{3.12}$$

*and* $\qquad \|\sigma_{\lambda\alpha}(s) - x\|_p \leq \dfrac{\|(B(x)+\omega(s))-A(x)\|_p}{\alpha - \|B\|_{op}}$, *for any* $s \in W_{\lambda\alpha}$ *and for* $x \in l_p$. $\qquad (3.13)$

*Proof.* For the given positive numbers $p$ and $q$ with $\dfrac{1}{p} + \dfrac{1}{q} = 1$, both $l_p$ and $l_q$ are uniformly convex and uniformly smooth Banach spaces. Hence, they are Asplund Banach spaces.

In Corollary 2.2 of the AMZ Theorem, let $X = Y = l_p$, which is Asplund. Further, let $F = A$ with $A: l_p \to l_p$ and let $G(\cdot, \cdot): l_p \times S \to l_p$, which is defined by $G(x, s) = B(x) + \omega(s)$, for $(x, s) \in l_p \times S$, being single-valued mappings. Since $A: l_p \to l_p$ is a linear and continuous single-valued mapping, then $A$ satisfies condition (A1) in the Corollary 2.2. By the condition (a1) in this theorem, and by part (iii) of Theorem 3.2, for every $x \in l_p$, we have $\hat{\alpha}(A, x, A(x)) > 0$, which satisfies

$$\hat{\alpha}(A, x, A(x)) = \min\{\|A^T(y)\|_q : y \in l_q, \|y\|_q = 1\},$$

and $\qquad 0 < \|B\|_{op} < \hat{\alpha}(A, x, A(x)) \leq \|A^T\|_{op} \leq 1$, for any $x \in l_p$.

In Corollary 2.2, for every $x \in l_p$, we take neighborhoods $U \subset l_p$ of $x$, $V \subset l_p$ of $A(x)$ to be $l_p$. That is, let $U = V = l_p$. For any fixed $s \in S$, by the continuity of $B$, we have

$$\|G(u, s) - G(v, s)\|_p = \|(B(u) + \omega(s)) - (B(v) + \omega(s))\|_p$$

$$= \|B(u) - B(v)\|_p \leq \|B\|_{op} \|u - v\|_p, \text{ for any } u, v \in l_p.$$

This implies that, for any $s \in S$, the function $G(\cdot, s): l_p \to l_p$ is Lipschitz-like on $U = l_p$ relative to $V = l_p$ with the uniform modulus $\|B\|_{op}$. By condition (a1) in this theorem, it satisfies

$$0 < \|B\|_{op} < \hat{\alpha}(A, x, A(x)) \leq 1, \text{ for any } x \in l_p.$$

By condition (a2), the function $s \to \omega(s)$ is lower semicontinuous at $\bar{s}$. This implies that, for any given $x \in l_p$, the function $G(x, s) = B(x) + \omega(s): l_p \times S \to l_p$ is lower semicontinuous at $\bar{s}$. Hence, condition (A2) in Corollary 2.2 is satisfied. By condition (a2) again, the Lipschitzian modulus $\|B\|_{op}$ of $G(\cdot, s)$ satisfies $0 < \|B\|_{op} < \hat{\alpha}(F, \bar{x}, F(\bar{y}))$. This implies that condition (A3) in Corollary 2.2 is satisfied. Hence, for any $\lambda, \alpha$ with $\hat{\alpha}(A, x, A(x)) \geq \lambda > \alpha > \|B\|_{op}$, there exist a neighborhood $W_{\lambda\alpha} \subset S$ of $\bar{s}$ and a single-valued mapping $\sigma_{\lambda\alpha}: W_{\lambda\alpha} \to l_p$ such that, (3.12) is satisfied. That is,

$$A(\sigma_{\lambda\alpha}(s)) = B(\sigma_{\lambda\alpha}(s)) + \omega(s), \text{ for every } s \in W_{\lambda\alpha},$$

and $\qquad \|\sigma_{\lambda\alpha}(s) - x\|_p \leq \dfrac{\text{dist}(A(x), B(x)+\omega(s))}{\alpha - \|B\|_{op}}$, for any $s \in W_{\lambda\alpha}$ and for $x \in l_p$.

Since $\text{dist}(A(x), B(x) + \omega(s)) = \|(B(x) + \omega(s)) - A(x)\|_p$, then (3.13) is proved. This theorem is completely proved. □

Let $C = (c_{ij})_{i,j=1}^{\infty}$ be a real $\infty \times \infty$ diagonal matrix. We need the following notations.

$$m_C = \inf\{|c_{ii}|: i = 1, 2, \dots\} \quad \text{and} \quad M_C = \sup\{|c_{ii}|: i = 1, 2, \dots\}.$$

**Lemma 3.5**. *Let $A = (a_{ij})_{i,j=1}^{\infty}$ be a real $\infty \times \infty$ diagonal matrices. If $M_A < \infty$, then,*

  (i) *The mapping $A$ is linear and continuous with $\|A\|_{\text{op}} = M_A$;*
  (ii) *$A$ is Fréchet differentiable at every point in $l_p$ and $\nabla(A)(x) = A$, for any $x \in l_p$;*
  (iii) *The Mordukhovich derivative of $A$ satisfies that $\widehat{D}^*(A)(x, A(x)) = A^T$, for any $x \in l_p$;*
  (iv) *Suppose that $M_A \le 1$, then the covering constant for $A$ is constant in $l_p$ with*

$$\hat{\alpha}(A, x, A(x)) = m_A, \text{ for any } x \in l_p.$$

*Proof.* We only sketch the proof of (iv). By $\inf\{|a_{ii}|: i = 1, 2, \dots\} < \infty$, there is a subsequence $\{a_{i_k i_k}\}_{k=1}^{\infty}$ of $\{a_{ii}\}_{i=1}^{\infty}$ such that $\inf\{|a_{i_k i_k}|: k = 1, 2, \dots\} = m_A$, as $k \to \infty$. Let $y^{(i_k)} \in l_p$ be defined as in Theorem 3.4, which satisfies that

$$\|A^T(y^{(i_k)})\|_q = \|y^{(i_k)} A^T\|_q = |a_{i_k i_k}|, \text{ for } k = 1, 2, \dots,$$

and
$$\inf\{\|A^T(y^{(i_k)})\|_q : k = 1, 2, \dots\} = \inf\{\|A^T(y)\|_q : y \in l_q, \|y\|_q = 1\}.$$

Notice that $\widehat{D}^*(A)(x, A(x)) = A^T$, which is independent from $x \in l_p$. This implies that

$$\hat{\alpha}(A, x, A(x)) = \inf\{\|A^T(y)\|_q : y \in l_q, \|y\|_q = 1\}$$

$$= \inf\{\|A^T(y^{(i_k)})\|_q : k = 1, 2, \dots\} = \inf\{|a_{i_k i_k}| : k = 1, 2, \dots\} = m_A, \text{ for any } x \in l_p. \quad □$$

Notice that the condition (iv) that $M_A \le 1$ is only for easier calculation of $\hat{\alpha}(A, x, A(x))$, in which

$$x \in \mathbb{B}_X(\bar{x}, \eta) \implies A(x) \in \mathbb{B}_X(A(\bar{y}), \eta).$$

In case, since (iv) also follows from the covering criterion, then the condition $M_A \le 1$ is excessive.

**Proposition 3.6**. *Let $(S, \tau, \mu)$ be a topological probability space. Let $A = (a_{ij})_{i,j=1}^{\infty}$ and $B = (b_{ij})_{i,j=1}^{\infty}$ be real $\infty \times \infty$ diagonal matrices. Let $\omega: S \to l_p$ be a $\tau$-measurable single-valued mapping. Suppose that $A, B$ and $\omega$ satisfy the following conditions.*

  (a1) $0 < M_B < m_A \le M_A \le 1$;
  (a2) *The function $s \to \omega(s)$ is lower semicontinuous at $\bar{s}$.*

*Then, for any $\lambda, \alpha$ with $M_B < \alpha < \lambda \le m_A$, there exist a neighborhood $W_{\lambda\alpha} \subset S$ of $\bar{s}$ and a single-valued mapping $\sigma_{\lambda\alpha}: W_{\lambda\alpha} \to l_p$ such that*

$$A(\sigma_{\lambda\alpha}(s)) = B(\sigma_{\lambda\alpha}(s)) + \omega(s), \text{ for every } s \in W_{\lambda\alpha}, \tag{3.14}$$

and
$$\|\sigma_{\lambda\alpha}(s) - x\|_p \leq \frac{\|(B(x)+\omega(s))-A(x)\|_p}{\alpha - \|B\|_{op}}, \text{ for any } s \in W_{\lambda\alpha} \text{ and for any } x \in l_p. \quad (3.15)$$

*Proof.* Notice $A^* = A^T = A$. By Lemma 3.5, we have $\hat{\alpha}(A, x, A(x)) = m_A$, $M_A = \|A\|_{op}$ and $M_B = \|B\|_{op}$, the proof of this proposition is similar to the proof of Theorem 3.4 and it is omitted here. □

As a matter of fact, for diagonal matrices $A$ and $B$, the conditions in Proposition 3.6 are sufficient condition to ensure the existence of solutions to the locally or globally stochastic infinitely dimensional systems of linear equations with respect to $A$, $B$ and $\omega$ (That is (3.14), the first part of the results of Proposition 3.6). They are not necessary conditions. We have the following counterexample (only for satisfying (3.14) and not for satisfying (3.15)).

**Example 3.7.** Let $A = (a_{ij})_{i,j=1}^{\infty}$ and $B = (b_{ij})_{i,j=1}^{\infty}$ be real $\infty \times \infty$ diagonal matrices. Let $\omega: S \to l_p$ be a single-valued mapping. Let $\omega(s) = (\omega_1(s), \omega_2(s), \ldots)$, for any $s \in S$, in which $\omega_k: S \to \mathbb{R}$ is a real valued function, for $k = 1, 2, \ldots$. Suppose that $A$ and $B$ satisfy the conditions:

$$0 < \inf\{|a_{ii}|: i = 1, 2, \ldots\} \leq \sup\{|a_{ii}|: i = 1, 2, \ldots\} < \inf\{|b_{ii}|: i = 1, 2, \ldots\} \leq \sup\{|b_{ii}|: i = 1, 2, \ldots\} \leq 1.$$

That is that $0 < \hat{\alpha}(A, x, A(x)) \leq \|A\|_{op} < \|B\|_{inf} \leq \|B\|_{op} \leq 1$. This implies that $A$, $B$ do not satisfy all of the conditions in Theorem 3.4. However, $A - B$ satisfies $|b_{ii} - a_{ii}| \geq \|B\|_{min} - \|A\|_{op} > 0$, for $i = 1, 2, \ldots, n$. Hence the linear and continuous mapping $A - B$ is a pointwise multiplication operator on $\mathbb{R}^n$. For $\omega(s) = (\omega_1(s), \omega_2(s), \ldots)$, define $\sigma(s) = (\sigma_1(s), \sigma_2(s), \ldots)$ by

$$\sigma_i(s) = \frac{1}{a_{ii} - b_{ii}} \omega_i(s), \text{ for } i = 1, 2, \ldots, \text{ for every } s \in S.$$

That is $\sigma(s) = (A - B)^{-1}\omega(s)$, for every $s \in S$, which implies that $A(\sigma(s)) = B(\sigma(s)) + \omega(s)$, for every $s \in S$. Hence, $\sigma$ satisfies (3.14), and $\sigma$ is a solution to the globally stochastic system of linear equations with respect to $A$, $B$ and $\omega$. However, with $\|B\|_{op}$ being the uniform modulus $\beta$ of the mapping $B(x) + \omega(s)$, which does not satisfy condition (a1) in Proposition 3.6, then, we show that the above solution $\sigma$ does not satisfy the inequality (3.15) in Proposition 3.6. Since $\hat{\alpha}(A, x, A(x)) > 0$, we have

$$\|A^{-1}\|_{op} = \frac{1}{\|A\|_{min}} = \frac{1}{\hat{\alpha}(A,x,A(x))}.$$

For any $x \in l_p$, by $\hat{\alpha}(A, x, A(x)) < \|B\|_{op}$ (condition (a1) in Proposition 3.6 is not satisfied), we obtain

$$\|\sigma(s) - x\|_p = \|(\sigma(s)B + \omega(s) - xA)A^{-1}\|_p \nleq \frac{\|A(x)-(B(x)+\omega(s))\|_p}{\hat{\alpha}(A,x,A(x)) - \|B\|_{op}}, \text{ for any } s \in S.$$

This implies that (3.15) does not hold.

## 4. Stochastic Integral-Equations

### 4.1. Integral Operator with Kernels

Let $a$ and $b$ be real numbers or infinity with $-\infty \leq a < b \leq \infty$. In case if $a = -\infty$ and $b = \infty$, then we consider $[a, b]$ to be $(-\infty, \infty)$. We write $[a, b] \times [a, b] = [a, b]^2$. Let $(L_2[a, b], \|\cdot\|_2)$ and $(L_2[a, b]^2, \|\cdot\|)$ denote the real Hilbert spaces of square integrable real valued functions defined on $[a, b]$ and $[a, b]^2$, respectively. Their origins are denoted by $\theta$. Let $k \in L_2[a, b]^2$ satisfying

$$\|k\| = \left(\int_{[a,b]^2} |k(u,v)|^2 \, dv \, du\right)^{\frac{1}{2}} < \infty. \tag{4.1}$$

Then, $k$ induces a mapping $K: L_2[a,b] \to L_2[a,b]$, which is defined, for any $f \in L_2[a,b]$, by

$$K(f)(u) = \int_a^b k(u,v) f(v) \, dv, \text{ for any } u \in [a,b]. \tag{4.2}$$

This mapping $K$ is an integral operator corresponding to (or induced by) the kernel $k \in L_2[a,b]^2$ and $K$ is a linear and continuous mapping from $L_2[a,b]$ to itself. Let $\|K\|_{\text{op}}$ denote its operator norm satisfying

$$\|K\|_{\text{op}} = \sup\{\|K(f)\|_2 : \|f\|_2 = 1\} = \sup\left\{\left(\int_a^b \left|\int_a^b k(u,v) f(v) \, dv\right|^2 du\right)^{\frac{1}{2}} : \int_a^b |f(v)|^2 \, dv = 1\right\}$$

$$\leq \left(\int_a^b \int_a^b |k(u,v)|^2 \, dv \, du\right)^{\frac{1}{2}} = \|k\| < \infty. \tag{4.3}$$

Let $(S, \tau, \mu)$ be a topological probability space as defined in the previous section. Let $\lambda$ be a real number. Let $\omega: S \to L_2[a,b]$ be a $\tau$-measurable single-valued mapping (that is considered as a noise). Let $\bar{s} \in S$. If there is a neighborhood (an event) $W \subset S$ of $\bar{s}$ and a single-valued mapping $\sigma: W \to L_2[a,b]$ such that

$$\lambda(\sigma(s))(u) = \int_a^b k(u,v)(\sigma(s))(v) \, dv + \omega(s)(u), \text{ for any } s \in W \text{ and for any } u \in [a,b], \tag{4.4}$$

then $\sigma$ is called a solution to the *locally stochastic integral-equations on $W$* around the possible outcome $\bar{s}$ and with respect to the kernel $k \in L_2[a,b]^2$ and the noise $\omega$. In particular, in the above inclusion property, if the neighborhood of $\bar{s}$ is the whole space $S$; that is, if

$$\lambda(\sigma(s))(u) = \int_a^b k(u,v)(\sigma(s))(v) \, dv + \omega(s)(u), \text{ for any } s \in S \text{ and for any } u \in [a,b], \tag{4.5}$$

then $\sigma$ is called a solution to the *globally stochastic integral-equations on $S$* around the possible outcome $\bar{s}$ and with respect to the kernel $k \in L_2[a,b]^2$ and the noise $\omega$.

In the proof of the existence of solutions for some locally or globally stochastic integral-equations, we need the results of Proposition 6.1 in [15]. We review it below.

**Proposition 6.1 in [15].** *Let $(X, \|\cdot\|)$ be a real Banach space with dual space $(X^*, \|\cdot\|_*)$. Let $I_X$ be the identity mapping in $X$. For any real number $\lambda$, the linear and continuous mapping $\lambda I_X: X \to X$ satisfies*

(i) *$\lambda I_X$ is Fréchet differentiable at every point in $X$ such that $\nabla(\lambda I_X)(x) = \lambda I_X$, for any $x \in X$;*

(ii) *The Mordukhovich derivative of $\lambda I_X$ satisfies that $\widehat{D}^*(\lambda I_X)(x, \lambda x) = \lambda I_{X^*}$, for any $x \in X$.*

(iii) *In addition, if $|\lambda| \leq 1$, then the covering constant for $\lambda I_X$ is constant in $X$ with*

$$\widehat{\alpha}(\lambda I_X, x, \lambda x) = |\lambda|, \text{ for any } x \in X.$$

**Theorem 4.1.** *Let $(S, \tau, \mu)$ be a topological probability space. Let $k \in L_2[a,b]^2$ and let $\omega: S \to L_2[a,b]$ be a single-valued mapping. Let $\lambda$ be a real number and $\bar{s} \in S$. Suppose that $k, \lambda$ and $\omega$ satisfy the following conditions.*

(a1) $0 < \|k\| < |\lambda| \leq 1$;

(a2) *The function $s \to \omega(s)$ is lower semicontinuous at $\bar{s}$.*

*Then, for any $\alpha$ with $\|k\| < \alpha < |\lambda|$, there exist a neighborhood $W_\alpha \subset S$ of $\bar{s}$ and a single-valued mapping $\sigma_\alpha : W_\alpha \to L_2[a,b]$ such that*

$$\lambda(\sigma_\alpha(s))(u) = \int_a^b k(u,v)(\sigma_\alpha(s))(v)dv + \omega(s)(u), \text{ for a.a } s \in W_\alpha \text{ and for a. a. } u \in [a,b], \quad (4.6)$$

*and* $\quad \|\sigma_\alpha(s) - f\|_2 \leq \dfrac{\left(\int_a^b \left|\int_a^b k(u,v)f(v)dv + \omega(s)(u) - \lambda f(u)\right|^2 du\right)^{\frac{1}{2}}}{\alpha - \|k\|}, \text{ for any } s \in W_\alpha \text{ and any } f \in L_2[a,b].$ (4.7)

*Proof.* Let $I$ denote the identity mapping in $L_2[a,b]$. For any real number $\lambda$ with $|\lambda| \leq 1$, by Proposition 6.1 in [15], we have that the covering constant for $\lambda I$ is constant in $L_2[a,b]$ with

$$\hat{\alpha}(\lambda I, f, \lambda f) = |\lambda|, \text{ for any } f \in L_2[a,b]. \quad (4.8)$$

In Corollary 2.2, let $F: L_2[a,b] \to L_2[a,b]$ be defined by $F = \lambda I$. It is clear that $F$ satisfies condition (A1) in Corollary 2.2. Let $K: L_2[a,b] \to L_2[a,b]$ be the linear and continuous integral operator induced by (or corresponding to) the kernel $k \in L_2[a,b]^2$. By condition (a1), $K$ satisfies $0 < \|K\|_{op} \leq \|k\| < |\lambda| \leq 1$. In Corollary 2.2, define $G: L_2[a,b] \times S \to L_2[a,b]$ by $G(h,s) = K(h) + \omega(s)$, for $(h,s) \in L_2[a,b] \times S$. Then, for any fixed $s \in S$, the function $G(\cdot, s): L_2[a,b] \to L_2[a,b]$ is Lipschitz-like on the whole space $L_2[a,b]$ relative to the whole space $L_2[a,b]$ with the uniform modulus $\|K\|_{op}$, which is constant with respect to $(h,s) \in L_2[a,b] \times S$. By condition (a1) in this theorem and by (4.8), it satisfies that

$$0 < \|K\|_{op} \leq \|k\| < |\lambda| = \hat{\alpha}(\lambda I, f, \lambda f) \leq 1, \text{ for any } f \in L_2[a,b].$$

By condition (a2) in this theorem, the function $s \to \omega(s)$ is lower semicontinuous at $\bar{s}$. This implies that, for any given $h \in L_2[a,b]$, the function $G(h,s) = K(h) + \omega(s): S \to L_2[a,b]$ is lower semicontinuous at $\bar{s}$. Hence, conditions (A2) and (A3) in Corollary 2.2 are satisfied. As a Hilbert space, $L_2[a,b]$ is an Asplund Banach space. Then, by Corollary 2.2, for any $\alpha$ with $\hat{\alpha}(A,x,A(x)) = |\lambda| > \alpha > \|k\| \geq \|K\|_{op}$, there exists a neighborhood $W_\alpha \subset S$ of $\bar{s}$ and a single-valued mapping $\sigma_\alpha: W_\alpha \to L_2[a,b]$ such that

$$\lambda(\sigma_\alpha(s))(u) = \int_a^b k(u,v)(\sigma_\alpha(s))(v)dv + \omega(s)(u), \text{ for any } s \in W_\alpha \text{ and for any } u \in [a,b].$$

Next, we check that $\sigma_\alpha(s)$ satisfies the inequality (4.7). To this end, we rewrite (4.6) as

$$(\sigma_\alpha(s))(u) = \tfrac{1}{\lambda} K(\sigma_\alpha(s))(u) + \tfrac{1}{\lambda}\omega(s)(u), \text{ for any } s \in W_\alpha \text{ and for any } u \in [a,b],$$

Then, by (4.6), for any $f \in L_2[a,b]$ (with $(\lambda I)(f) = \lambda f$), we have

$$\|\sigma_\alpha(s) - f\|_2 = \left\|\tfrac{1}{\lambda}K(\sigma_\alpha(s)) + \tfrac{1}{\lambda}\omega(s) - f\right\|_2 = \tfrac{1}{\lambda}\|K(\sigma_\alpha(s)) + \omega(s) - \lambda f\|_2$$

$$= \tfrac{1}{\lambda}\|K(\sigma_\alpha(s)) + \omega(s) - (\lambda I)(f)\|_2 \leq \dfrac{\|K(\sigma_\alpha(s)) - K(f)\|_2 + \|K(f) + \omega(s) - (\lambda I)(f)\|_2}{\lambda}$$

$$\leq \dfrac{\|K\|_{op}\|\sigma_\alpha(s) - f\|_2 + \|K(f) + \omega(s) - (\lambda I)(f)\|_2}{\lambda}.$$

This implies that

$$\|\sigma_\alpha(s) - f\|_2 (1 - \frac{\|K\|_{\text{op}}}{\lambda}) \leq \frac{\|K(f) + \omega(s) - (\lambda I)(f)\|_2}{\lambda}.$$

By $0 < \|K\|_{\text{op}} \leq \|k\| < \hat{\alpha}(\lambda I, f, \lambda f) \leq 1$, and $\|K\|_{\text{op}} \leq \|k\| < \alpha < \lambda = \hat{\alpha}(\lambda I, f, \lambda f)$, we obtain

$$\|\sigma_\alpha(s) - f\|_2 \leq \frac{\|K(f) + \omega(s) - (\lambda I)(f)\|_2}{\lambda - \|K\|_{\text{op}}} \leq \frac{\|(K(f) + \omega(s)) - \lambda f\|_2}{\alpha - \|K\|_{\text{op}}}$$

$$= \frac{\left(\int_a^b \left|\int_a^b k(u,v)f(v)dv + \omega(s)(u) - \lambda f(u)\right|^2 du\right)^{\frac{1}{2}}}{\alpha - \|K\|_{\text{op}}} \leq \frac{\left(\int_a^b \left|\int_a^b k(u,v)f(v)dv + \omega(s)(u) - \lambda f(u)\right|^2 du\right)^{\frac{1}{2}}}{\alpha - \|k\|},$$

which proves (4.7). □

In particular, in Theorem 4.1, let $\omega: S \to L_2[a,b]$ be the $\theta$-functional. Then, we obtain the existence of solutions for some integral equations by the following corollary.

**Corollary 4.2.** *Let $(S, \tau, \mu)$ be a topological probability space. Let $k \in L_2[a,b]^2$ and let $\omega: S \to L_2[a,b]$ be the $\theta$-functional. Let $\lambda$ be a real number. Suppose that $k$ and $\lambda$ satisfy the following conditions.*

(a1) $0 < \|k\| < |\lambda| \leq 1$.

*Then, for any $\bar{s} \in S$ and for any $\alpha$ with $1 \geq |\lambda| > \alpha > \|k\|$, there exist a neighborhood $W_\alpha \subset S$ of $\bar{s}$ and a single-valued mapping $\sigma_\alpha: W_\alpha \to L_2[a,b]$ such that*

$$\lambda(\sigma_\alpha(s))(u) = \int_a^b k(u,v)(\sigma_\alpha(s))(v)dv, \text{ for any } s \in W_\alpha \text{ and for any } u \in [a,b],$$

*and*
$$\|\sigma_\alpha(s) - f\|_2 \leq \frac{\left(\int_a^b \left|\int_a^b k(u,v)f(v)dv - \lambda f(u)\right|^2 du\right)^{\frac{1}{2}}}{\alpha - \|k\|}, \text{ for any } s \in W_\alpha \text{ and any } f \in L_2[a,b].$$

Another special case of Theorem 4.1 is either both the kernel function $k$ and the noise $\omega$ take the especially simple form of tensor products, or, at least the kernel function $k$ is a tensor product of two single variable functions.

**Corollary 4.3.** *Let $(S, \tau, \mu)$ be a topological probability space. Let $h, g \in L_2[a,b]$ and let $\omega: S \to L_2[a,b]$ be a single-valued mapping. Let $\lambda$ be a real number. Suppose that $h, g, \omega$ and $\lambda$ satisfy the following conditions.*

(a1) $0 < \|g\|_2 \|h\|_2 < |\lambda| \leq 1$;

(a2) *The function $s \to \omega(s)$ is lower semicontinuous at $\bar{s}$.*

*Then, for this $\bar{s} \in S$ and for any $\alpha$ with $|\lambda| > \alpha > \|g\|_2 \|h\|_2$, there exist a neighborhood $W_\alpha \subset S$ of $\bar{s}$ and a single-valued mapping $\sigma_\alpha: W_\alpha \to L_2[a,b]$ such that*

$$\lambda(\sigma_\alpha(s))(u) = g(u) \int_a^b h(v)(\sigma_\alpha(s))(v)dv + \omega(s)(u), \text{ for any } s \in W_\alpha \text{ and for any } u \in [a,b], \quad (4.9)$$

*and* $\|\sigma_\alpha(s) - f\|_2 \leq \dfrac{\left(\int_a^b \left|g(u)\int_a^b h(v)f(v)dv + \omega(s)(u) - \lambda f(u)\right|^2 du\right)^{\frac{1}{2}}}{\alpha - \|g\|_2 \|h\|_2},$ *for any $s \in W_\alpha$ and $f \in L_2[a,b]$.* (4.10)

*Moreover, under condition* $1 \geq |\lambda| > \alpha > \|g\|_2 \|h\|_2$, $\sigma_\alpha(s)$ *in (4.9) has the following property*

$$\int_a^b h(v)(\sigma_\alpha(s))(v)dv = \frac{\int_a^b h(v)\omega(s)(v)dv}{\lambda - \int_a^b g(v)h(v)dv}, \text{ for any } s \in W_\alpha, \qquad (4.11)$$

*and* $\sigma_\alpha(s)$ *in (4.9) has the following precisely explicit representation*

$$(\sigma_\alpha(s))(u) = \frac{1}{\lambda} \frac{\int_a^b h(v)\omega(s)(v)dv}{\lambda - \int_a^b g(v)h(v)dv} g(u) + \frac{1}{\lambda}\omega(s)(u), \text{ for any } s \in W_\alpha \text{ and for any } u \in [a,b]. \qquad (4.12)$$

*Proof.* Define $k(u,v) = g(u)h(v)$, for any $(u,v) \in [a,b]^2$. $k$ is a tensor product kernel with $k \in L_2[a,b]^2$ and $\|k\| = \|g\|_2 \|h\|_2$. Then, (4.9) and (4.10) follow from Theorem 4.1 and Corollary 4.2 immediately. We only need to show (4.11). At first, by condition (a1), we have $\left|\int_a^b g(v)h(v)dv\right| \leq \|g\|_2\|h\|_2 < |\lambda|$. This implies that the fraction in (4.11) is well-defined. Both sides of (4.9) are multiplied by $h(v)$ and taken integral on $[a,b]$, for any $s \in W_\alpha$, we have

$$\lambda \int_a^b h(u)(\sigma_\alpha(s))(u)du = \int_a^b g(u)h(u)du \int_a^b h(v)(\sigma_\alpha(s))(v)dv + \int_a^b h(u)\omega(s)(u)du.$$

By (4.9) and (4.11), this proves (4.12). □

**Example 4.4.** Let $(S, \tau, \mu)$ be the topological probability space with $S = [0,1]$. Let $h, g \in L_2[-1,1]$ with

$$g(u) = u^2 \text{ and } h(v) = v^4, \text{ for any } u, v \in [-1,1].$$

Let $\omega(s)(v) = s^2 v^2$, for any $s \in [0,1]$ and $v \in [-1,1]$. We have $\|g\|_2 \|h\|_2 = \frac{2}{3\sqrt{5}}$. Define $k(u,v) = g(u)h(v) = u^2 v^4$, for any $(u,v) \in [a,b]^2$. $k$ is a tensor product kernel with $k \in L_2[a,b]^2$ and $\|k\| = \|g\|_2 \|h\|_2 = \frac{2}{3\sqrt{5}}$. Define a linear and continuous mapping $K$ on $L_2[-1,1]$ corresponding to the kernel $k$. Then, for any real number $\lambda$ with $1 \geq |\lambda| > \|g\|_2 \|h\|_2 = \frac{2}{3\sqrt{5}}$ and for any $s \in [0,1]$, by (4.9), we have a solution $\sigma_\alpha$ of the stochastic integral equation

$$\lambda(\sigma_\alpha(s))(u) = u^2 \int_{-1}^{1} v^4(\sigma_\alpha(s))(v)dv + s^2 u^2, \text{ for any } s \in [0,1] \text{ and for any } u \in [-1,1].$$

We calculate that $\int_{-1}^{1} g(v)h(v)dv = \int_{-1}^{1} v^6 dv = \frac{2}{7}$ and

$$\int_{-1}^{1} h(v)\omega(s)(v)dv = \int_{-1}^{1} v^4(\omega(s))(v)dv = \frac{2s^2}{7}, \text{ for any } s \in [0,1].$$

Substituting the above results into (4.12), the solution $\sigma_\alpha(s)$ has the following explicit representation.

$$(\sigma_\alpha(s))(u) = \frac{1}{\lambda} \frac{2s^2 u^2}{7\lambda - 2} + \frac{1}{\lambda} s^2 u^2, \text{ for any } s \in [0,1] \text{ and for any } u \in [-1,1].$$

**Example 4.5.** Let $(S, \tau, \mu)$ be the topological probability space with $S = [0,1]$ as used in Example 4.4. Let again $h, g \in L_2[-1,1]$ with $g(u) = u^2$ and $h(v) = v^4$, for any $u, v \in [-1,1]$. Let $\omega(s)(v) = s^2 \sin v$, for any $s \in [0,1]$ and $v \in [-1,1]$. We have $\|g\|_2 \|h\|_2 = \frac{2}{3\sqrt{5}}$. Let $k$ and $K$ be defined as in the last example. Then, for any real number $\lambda$ with $1 \geq |\lambda| > \|g\|_2 \|h\|_2 = \frac{2}{3\sqrt{5}}$ and for any $s \in [0,1]$, by (4.12), we have a solution $\sigma_\alpha$ of the stochastic integral equation with respect to $K$

$$\lambda(\sigma_\alpha(s))(u) = u^2 \int_{-1}^{1} v^4(\sigma_\alpha(s))(v)dv + s^2\sin u, \text{ for any } s \in [0, 1] \text{ and for any } u \in [-1,1].$$

We calculate that $\int_{-1}^{1} h(v)\omega(s)(v)dv = \int_{-1}^{1} v^4 s^2 \sin v\, dv = 0$, for any $s \in [0, 1]$. Substituting the above results into (4.12), the solution $\sigma_\alpha(s)$ has the following explicit representation.

$$(\sigma_\alpha(s))(u) = \frac{1}{\lambda} s^2 \sin u, \text{ for any } s \in [0, 1] \text{ and for any } u \in [-1,1].$$

**Example 4.6.** Let $(S, \tau, \mu)$ be the topological probability space with $S = [0, 1]$. Let $h, g \in L_2(\mathbb{R})$ with

$$g(u) = \sqrt{\frac{1}{4\pi}\frac{1}{1+u^2}} \text{ and } h(v) = \sqrt{\frac{1}{\sqrt{2\pi}}} e^{-\frac{v^2}{2}}, \text{ for any } u, v \in (-\infty, \infty).$$

Define $\omega(s)$, for any $s \in [0,1]$, by

$$\omega(s)(v) = \begin{cases} sv, & \text{if } |v| \le 1, \\ 0, & \text{if } |v| > 1, \end{cases} \text{ for any } v \in (-\infty, \infty).$$

We have $\|g\|_2 \|h\|_2 = \frac{1}{2}$. Let $k(u, v) = g(u)h(v)$, for any $(u, v) \in \mathbb{R}^2$. $k$ is a tensor product kernel with $k \in L_2(\mathbb{R}^2)$ and $\|k\| = \|g\|_2 \|h\|_2 = \frac{1}{2}$. Define the linear and continuous mapping $K$ on $L_2(\mathbb{R})$ induced by (corresponding to) the kernel $k$. Then, for any real number $\lambda$ with $1 \ge |\lambda| > \|g\|_2 \|h\|_2 = \frac{1}{2}$ and for any $s \in [0, 1]$, by (4.12), we have a solution $\sigma_\alpha$ of the stochastic integral equation with respect to $K$

$$\lambda(\sigma_\alpha(s))(u) = \begin{cases} u^2 \int_{-\infty}^{\infty} v^4(\sigma_\alpha(s))(v)dv + su, & \text{if } |u| \le 1, \\ u^2 \int_{-\infty}^{\infty} v^4(\sigma_\alpha(s))(v)dv, & \text{if } |u| > 1, \end{cases} \text{ for any } s \in [0, 1] \text{ and } u \in (-\infty, \infty).$$

We calculate

$$\int_{-\infty}^{\infty} h(v)\omega(s)(v)dv = \int_{-\infty}^{\infty} \sqrt{\frac{1}{\sqrt{2\pi}}} e^{-\frac{v^2}{2}} \omega(s)(v)dv = \int_{-1}^{1} \sqrt{\frac{1}{\sqrt{2\pi}}} e^{-\frac{v^2}{2}} s^2 v\, dv = 0, \text{ for any } s \in [0, 1].$$

Substituting the above results into (4.12), the solution $\sigma_\alpha(s)$ has the following explicit representation.

$$(\sigma_\alpha(s))(u) = \frac{1}{\lambda}\omega(s)(u) = \begin{cases} \frac{su}{\lambda}, & \text{if } |u| \le 1, \\ 0, & \text{if } |u| > 1, \end{cases} \text{ for any } s \in [0, 1] \text{ and for any } u \in (-\infty, \infty).$$

**4.2. Integral Operator with Orthonormal Basis**

Let $(L_2[a, b], \|\cdot\|_2)$ and $(L_2[a, b]^2, \|\cdot\|)$ be the Hilbert spaces used in the previous subsection. Suppose that $L_2[a, b]$ has an orthonormal (Schauder) basis $\{e_n\}_{n=1}^{\infty}$ that satisfies

$$\int_a^b e_m(v)e_n(v)dv = \begin{cases} 1, & \text{if } m = n, \\ 0, & \text{if } m \ne n, \end{cases}$$

For any $f \in L_2[a, b]$, let $\{f_n\}_{n=1}^{\infty}$ be the sequence of coefficients of $f$ with respect to the orthonormal Schauder basis $\{e_n\}_{n=1}^{\infty}$ of $L_2[a, b]$ such that

$$f_n = \langle f, e_n \rangle = \int_a^b f(v)e_n(v)dv, \text{ for } n = 1, 2, \ldots. \tag{4.13}$$

and $\quad f(v) = \sum_{n=1}^{\infty} f_n e_n(v)$, for any $v \in [a, b]$ satisfying $\|f\|_2 = (\sum_{n=1}^{\infty} f_n^2)^{\frac{1}{2}} < \infty.$ (4.14)

For any positive integers $m$ and $n$, define

$$e_{mn}(u, v) = (e_m \times e_n)(u, v) = e_m(u) e_n(v), \text{ for any } (u, v) \in [a, b]^2. \quad (4.15)$$

Then, it is known that $\{e_{mn}\}_{m,n=1}^{\infty}$ forms an orthonormal (Schauder) basis of $L_2[a, b]^2$. This implies that, for any $k \in L_2[a, b]^2$, there is a square summable double sequence of real numbers $\{k_{mn}\}_{m,n=1}^{\infty}$, such that

$$k(u, v) = \sum_{m=1}^{\infty} \sum_{n=1}^{\infty} k_{mn} e_{mn}(u, v) = \sum_{m=1}^{\infty} \sum_{n=1}^{\infty} k_{mn} e_m(u) e_n(v), \text{ for any } u, v \in [a, b]. \quad (4.16)$$

$\{k_{mn}\}_{m,n=1}^{\infty}$ is the sequence of coefficients of $k$ with respect to the orthonormal Schauder basis $\{e_{mn}\}_{m,n=1}^{\infty}$ of $L_2[a, b]^2$. By (4.1), the $L_2[a, b]^2$-norm of $k$ satisfies

$$\|k\| = \left( \int_{[a,b]^2} (\sum_{m=1}^{\infty} \sum_{n=1}^{\infty} k_{mn} e_m(u) e_n(v))^2 \, dudv \right)^{\frac{1}{2}} = (\sum_{m,n=1}^{\infty} k_{mn}^2)^{\frac{1}{2}} < \infty. \quad (4.17)$$

Let $K: L_2[a, b] \to L_2[a, b]$ be the linear and continuous integral operator corresponding to (or induced by) the kernel $k$. With the representation (4.13), $K$ is defined, for any $f \in L_2[a, b], u \in [a, b]$, by

$$K(f)(u) = \int_a^b k(u, v) f(v) dv = \int_a^b (\sum_{m=1}^{\infty} \sum_{n=1}^{\infty} k_{mn} e_{mn}(u, v)) f(v) dv$$

$$= \int_a^b (\sum_{m=1}^{\infty} \sum_{n=1}^{\infty} k_{mn} e_m(u) e_n(v)) f(v) dv = \int_a^b (\sum_{m=1}^{\infty} \sum_{n=1}^{\infty} k_{mn} e_m(u) e_n(v)) (\sum_{j=1}^{\infty} f_j e_j(v)) dv$$

$$= \sum_{m=1}^{\infty} e_m(u) \int_a^b (\sum_{n=1}^{\infty} k_{mn} e_n(v)) (\sum_{j=1}^{\infty} f_j e_j(v)) dv = \sum_{m=1}^{\infty} (\sum_{n=1}^{\infty} k_{mn} f_n) e_m(u). \quad (4.18)$$

This mapping $K$ is an integral operator corresponding to the kernel $k \in L_2[a, b]^2$ and $K$ is a linear and continuous mapping from $L_2[a, b]$ to itself. By (4.17), the operator norm $\|K\|_{op}$ of $K$ satisfies

$$\|K\|_{op} = \sup\{\|K(f)\|_2 : \|f\|_2 = 1\} = \sup\left\{ \left( \int_a^b \left( \int_a^b k(u, v) f(v) dv \right)^2 du \right)^{\frac{1}{2}} : \sum_{j=1}^{\infty} f_j^2 = 1 \right\}$$

$$= \sup\left\{ \left( \int_a^b (\sum_{m=1}^{\infty} (\sum_{n=1}^{\infty} k_{mn} f_n) e_m(u))^2 du \right)^{\frac{1}{2}} : \sum_{j=1}^{\infty} f_j^2 = 1 \right\}$$

$$= \sup\left\{ (\sum_{m=1}^{\infty} (\sum_{n=1}^{\infty} k_{mn} f_n)^2)^{\frac{1}{2}} : \sum_{j=1}^{\infty} f_j^2 = 1 \right\} \leq \sup\left\{ (\sum_{m=1}^{\infty} (\sum_{n=1}^{\infty} k_{mn}^2) (\sum_{n=1}^{\infty} f_n^2))^{\frac{1}{2}} : \sum_{j=1}^{\infty} f_j^2 = 1 \right\}$$

$$= (\sum_{m=1}^{\infty} \sum_{n=1}^{\infty} k_{mn}^2)^{\frac{1}{2}} = \|k\| < \infty.$$

Let $A = \{a_n\}_{n=1}^{\infty}$ be a bounded sequence of real numbers, which defines a pointwise multiplication operator on $L_2[a, b]$ with respect to the orthonormal basis $\{e_n\}_{n=1}^{\infty}$ of $L_2[a, b]$. $A$ is defined by the following explicit representation with respect to the orthonormal basis $\{e_n\}_{n=1}^{\infty}$ of $L_2[a, b]$,

$$(A(f))(v) = \sum_{n=1}^{\infty} a_n f_n e_n(v), \text{ for any } v \in [a, b]. \quad (4.19)$$

This operator is a linear and continuous mapping from $L_2[a, b]$ to itself. For any $f(v) = \sum_{n=1}^{\infty} f_n e_n(v) \in L_2[a, b]$, $A(f) \in L_2[a, b]$. Similar, to (4.14), we have $\|A(f)\|_2 = (\sum_{n=1}^{\infty} a_n^2 f_n^2)^{\frac{1}{2}} < \infty$. This implies that

$$\|A\|_{\inf} = \inf\{\,|a_n|: n = 1, 2, \ldots\} \leq \|A\|_{\mathrm{op}} = \sup\{\,|a_n|: n = 1, 2, \ldots\} < \infty.$$

**Lemma 4.7.** *Let $A = \{a_n\}_{n=1}^{\infty}$ be a bounded sequence of real numbers that defines a pointwise multiplication operator $A$ on $L_2[a,b]$ by (4.19). Then, for any point $f \in L_2[a,b]$, we have*

(i) *$A$ is Fréchet differentiable at a point $f$ and $\nabla(A)(f) = A$;*

(ii) *The Mordukhovich derivative of $A$ satisfies that $\widehat{D}^*(A)(f, A(f)) = A$;*

(iii) *Suppose that $\sup\{\,|a_n|: n = 1, 2, \ldots\} \leq 1$, then the covering constant for $A$ is constant with*

$$\hat{\alpha}(A, f, A(f)) = \|A\|_{\inf} = \inf\{\|A(g)\|_2 : g \in \mathbb{S}_2\} = \inf\{|a_n| : n = 1, 2, \ldots\},$$

*and* $\qquad 0 \leq \hat{\alpha}(A, f, A(f)) = \|A\|_{\inf} \leq \|A\|_{\mathrm{op}} = \sup\{|a_n| : n = 1, 2, \ldots\} \leq 1.$

*Proof.* This Lemma actually follows from Lemma 3.5. However, we give a direct proof here. (i). The mapping $A: L_2[a,b] \to L_2[a,b]$ defined by (4.19) is a linear and continuous single-valued mapping. We have

$$\lim_{h \to f} \frac{A(f) - A(h) - A(f-h)}{\|f - h\|_2} = \theta, \text{ for any given } f \in L_2[a,b].$$

This proves (i). Part (ii) is proved by part (i) and Theorem 1.38 in [17]. Now, by (ii), we prove part (iii). Since $A$ satisfies that $\sup\{\,|a_j|: j = 1, 2, \ldots\} \leq 1$, we have $\|A\|_{\mathrm{op}} \leq 1$. Then, for any $f \in L_2[a,b]$ and for any $\eta > 0$, we have that

$$h \in \mathbb{B}_2(f, \eta) \implies A(h) \in \mathbb{B}_2(A(f), \eta), \text{ for } h \in L_2[a,b].$$

Since $L_2[a,b]$ is a Hilbert space, by the definition (4.19), the adjoint operator $A^*$ is $A$. By (ii), we calculate the covering constant for $A$ at an arbitrarily given point $f \in L_2[a,b]$.

$$\hat{\alpha}(A, f, A(f)) = \sup_{\eta > 0} \inf\{\|w\|_2 : w \in \widehat{D}^*(A)(h, A(h))(g), h \in \mathbb{B}_2(f, \eta), A(h) \in \mathbb{B}_2(A(f), \eta), \|g\|_2 = 1\}$$

$$= \sup_{\eta > 0} \inf\{\|A^*(g)\|_2 : \{A^*(g)\} = \widehat{D}^*(A)(h, A(h))(g), h \in \mathbb{B}_2(f, \eta), A(h) \in \mathbb{B}_2(A(f), \eta), \|g\|_2 = 1\}$$

$$= \sup_{\eta > 0} \inf\{\|A(g)\|_2 : g \in L_2[a,b], \|g\|_2 = 1\}$$

$$= \inf\{\|A(g)\|_2 : g \in L_2[a,b], \|g\|_2 = 1\}$$

$$= \inf\{|a_n| : n = 1, 2, \ldots\}$$

$$= \|A\|_{\inf}. \qquad \square$$

**Theorem 4.8.** *Let $(S, \tau, \mu)$ be a topological probability space. Let $\{e_n\}_{n=1}^{\infty}$ be an orthonormal Schauder basis of $L_2[a,b]$. Let $k = \sum_{m=1}^{\infty} \sum_{n=1}^{\infty} k_{mn} e_m e_n \in L_2[a,b]^2$ with*

$$k(u, v) = \sum_{m=1}^{\infty} \sum_{n=1}^{\infty} k_{mn} e_m(u) e_n(v), \text{ for every } (u, v) \in [a, b]^2.$$

*Let $\omega: S \to L_2[a,b]$ be a single-valued mapping with sequence of coefficients $\{(\omega(s))_n\}_{n=1}^{\infty}$, for any $s \in S$. Let $A = \{a_n\}_{n=1}^{\infty}$ be a bounded sequence of real numbers that induces a pointwise multiplication* operator A. Let $\bar{s} \in S$. Suppose that $k, A$ and $\omega$ satisfy the following conditions.

(a1) $0 < \|k\| < \|A\|_{\inf} \leq \|A\|_{op} \leq 1$;
(a2) The function $s \to \omega(s)$ is lower semicontinuous at $\bar{s}$.

Then, for any $\alpha$ with $\|A\|_{\inf} > \alpha > \|k\|$, there exist a neighborhood $W_\alpha \subset S$ of $\bar{s}$ and a single-valued mapping $\sigma_\alpha: W_\alpha \to L_2[a,b]$ with sequence of coefficients $\{(\sigma_\alpha(s))_n\}_{n=1}^\infty$, for any $s \in S$ such that

$$a_m(\sigma_\alpha(s))_m = \sum_{n=1}^\infty k_{mn}(\sigma_\alpha(s))_n + (\omega(s))_m, \text{ for any } s \in W_\alpha \text{ and for } m = 1, 2, \ldots. \quad (4.20)$$

And, for any $f \in L_2[a,b]$, we have

$$\sqrt{\sum_{m=1}^\infty ((\sigma_\alpha(s))_m - f_m)^2} \leq \frac{\sqrt{\sum_{m=1}^\infty (\sum_{n=1}^\infty k_{mn}f_n + (\omega(s))_m - a_m f_m)^2}}{\alpha - \|k\|}, \text{ for any } s \in W_\alpha. \quad (4.21)$$

*Proof.* Since $\omega: S \to L_2[a,b]$ is a single-valued mapping, then, for every $s \in S$, $\omega(s)$ is represented by its sequence of coefficients with respect to the orthonormal Schauder basis $\{e_n\}_{n=1}^\infty$ of $L_2[a,b]$.

$$\omega(s)(u) = \sum_{n=1}^\infty (\omega(s))_n e_n(u), \text{ for any } u \in [a,b]. \quad (4.22)$$

Let $A: L_2[a,b] \to L_2[a,b]$ be the pointwise multiplication operator corresponding to the bounded sequence of real numbers $\{a_n\}_{n=1}^\infty$, which is a linear and continuous mapping on $L_2[a,b]$. By Lemma 4.7 and the condition (a1) in this theorem, the covering constant for $A$ constant satisfies

$$0 < \hat{\alpha}(A, f, A(f)) = \|A\|_{\inf} \leq \|A\|_{op} \leq 1, \text{ for any } f \in L_2[a,b].$$

Let $K: L_2[a,b] \to L_2[a,b]$ be the linear and continuous integral operator corresponding to the kernel $k$. Define a single-valued mapping $G(\cdot, \cdot): L_2[a,b] \times S \to L_2[a,b]$ by

$$G(f, s) = K(f) + \omega(s), \text{ for any } (f, s) \in L_2[a,b] \times S.$$

Then, it is easy to check that, for any $s \in S$, the single-valued mapping $G(\cdot, s): L_2[a,b] \to L_2[a,b]$ is a Lipschitz mapping on the whole space $L_2[a,b]$ with the uniform modulus $\beta = \|K\|_{op} \leq \|k\|$. Hence, all conditions in Theorem 4.1 (or in Corollary 2.2) are satisfied for the considered single-valued mappings $A$, $G$ and the noise $\omega$. Then, by Theorem 4.1 (or by Corollary 2.2), for any $\alpha$ with $\|A\|_{\inf} > \alpha > \|k\|$, there exists a neighborhood $W_\alpha \subset S$ of $\bar{s}$ and a single-valued mapping $\sigma_\alpha: W_\alpha \to L_2[a,b]$ such that, for $s \in W_\alpha$,

$$A(\sigma_\alpha(s))(u) = \int_a^b k(u,v)(\sigma_\alpha(s))(v)dv + \omega(s)(u), \text{ for any } u \in [a,b]. \quad (4.23)$$

Similarly, to (4.22), for any $s \in W_\alpha$, we have

$$(\sigma_\alpha(s))(u) = \sum_{n=1}^\infty (\sigma_\alpha(s))_n e_n(u), \text{ for any } u \in [a,b]. \quad (4.24)$$

Notice that the linear and continuous mapping $A$ is a pointwise multiplication operator on $L_2[a,b]$ with respect to the orthonormal Schauder basis $\{e_n\}_{n=1}^\infty$ of $L_2[a,b]$. Substituting (4.22) and (4.24) into (4.23), we obtain

$$\sum_{m=1}^\infty a_m(\sigma_\alpha(s))_m e_m(u) = A(\sum_{m=1}^\infty (\sigma_\alpha(s))_m e_m(u)) = A(\sigma_\alpha(s))(u)$$

$$= \int_a^b k(u,v)(\sigma_\alpha(s))(v)dv + \omega(s)(u)$$

$$= \int_a^b \sum_{m=1}^\infty \sum_{n=1}^\infty k_{mn} e_{mn}(u,v) \sum_{j=1}^\infty (\sigma_\alpha(s))_j e_j(v)\, dv + \sum_{m=1}^\infty (\omega(s))_m e_m(u)$$

$$= \int_a^b \sum_{m=1}^\infty \sum_{n=1}^\infty k_{mn} e_m(u) e_n(v) \sum_{j=1}^\infty (\sigma_\alpha(s))_j e_j(v)\, dv + \sum_{m=1}^\infty (\omega(s))_m e_m(u)$$

$$= \sum_{m=1}^\infty e_m(u) \int_a^b \sum_{n=1}^\infty k_{mn} e_n(v) \sum_{j=1}^\infty (\sigma_\alpha(s))_j e_j(v)\, dv + \sum_{m=1}^\infty (\omega(s))_m e_m(u)$$

$$= \sum_{m=1}^\infty \left(\sum_{n=1}^\infty k_{mn} (\sigma_\alpha(s))_n\right) e_m(u) + \sum_{m=1}^\infty (\omega(s))_m e_m(u)$$

$$= \sum_{m=1}^\infty \left(\sum_{n=1}^\infty k_{mn} (\sigma_\alpha(s))_n + (\omega(s))_m\right) e_m(u), \text{ for any } u \in [a,b].$$

This implies that $a_m (\sigma_\alpha(s))_m = \sum_{n=1}^\infty k_{mn} (\sigma_\alpha(s))_n + (\omega(s))_m$, for $m = 1, 2, \ldots$. Hence, (4.20) is proved. Moreover, further, for any $f(v) = \sum_{n=1}^\infty f_n e_n(v) \in L_2[a,b]$, by (4.23), we have

$$\|\sigma_\alpha(s) - f\|_2 = \|A^{-1}(K(\sigma_\alpha(s)) + \omega(s)) - f\|_2$$

$$= \|A^{-1}(K(\sigma_\alpha(s)) + \omega(s) - A(f))\|_2 \leq \|A^{-1}\|_{op} \|K(\sigma_\alpha(s)) + \omega(s) - A(f)\|_2$$

$$\leq \|A^{-1}\|_{op} \|K(\sigma_\alpha(s)) - K(f)\|_2 + \|A^{-1}\|_{op} \|K(f) + \omega(s) - A(f)\|_2$$

$$\leq \|A^{-1}\|_{op} \|K\|_{op} \|\sigma_\alpha(s) - f\|_2 + \|A^{-1}\|_{op} \|K(f) + \omega(s) - A(f)\|_2.$$

This implies that $\|\sigma_\alpha(s) - f\|_2 (1 - \|A^{-1}\|_{op} \|K\|_{op}) \leq \|A^{-1}\|_{op} \|K(f) + \omega(s) - A(f)\|_2$. Since $\|A^{-1}\|_{op} = \frac{1}{\|A\|_{min}}$. We obtain $\|\sigma_\alpha(s) - f\|_2 \left(1 - \frac{\|K\|_{op}}{\|A\|_{inf}}\right) \leq \frac{1}{\|A\|_{inf}} \|K(f) + \omega(s) - A(f)\|_2$. This implies

$$\|\sigma_\alpha(s) - f\|_2 \leq \frac{\|K(f) + \omega(s) - A(f)\|_2}{\|A\|_{inf} - \|K\|_{op}}. \tag{4.25}$$

We calculate

$$\|\sigma_\alpha(s) - f\|_2 = \left\|\sum_{m=1}^\infty (\sigma_\alpha(s))_m e_m - \sum_{m=1}^\infty f_m e_m\right\|_2 = \sqrt{\sum_{m=1}^\infty ((\sigma_\alpha(s))_m - f_m)^2}, \tag{4.26}$$

and
$$(K(f))(u) = \int_a^b k(u,v) f(v)\, dv$$

$$= \int_a^b \sum_{m=1}^\infty \sum_{n=1}^\infty k_{mn} e_m(u) e_n(v) \sum_{j=1}^\infty f_j e_j(v)\, dv = \sum_{m=1}^\infty \left(\sum_{n=1}^\infty k_{mn} f_n\right) e_m(u).$$

This implies that

$$\|K(f) + \omega(s) - A(f)\|_2 = \left\|\sum_{m=1}^\infty \left(\sum_{n=1}^\infty k_{mn} f_n\right) e_m + \sum_{m=1}^\infty (\omega(s))_m e_m - \sum_{m=1}^\infty a_m f_m e_m\right\|_2$$

$$= \sqrt{\sum_{m=1}^\infty \left(\sum_{n=1}^\infty k_{mn} f_n + (\omega(s))_m - a_m f_m\right)^2}. \tag{4.27}$$

Substituting (4.26) and (4.27) into (4.25), for any $\alpha$ with $\|A\|_{inf} > \alpha > \|k\| \geq \|K\|_{op}$, we have

$$\sqrt{\sum_{m=1}^\infty ((\sigma_\alpha(s))_m - c_{f,m})^2} \leq \frac{\sqrt{\sum_{m=1}^\infty \left(\sum_{n=1}^\infty k_{mn} f_n + (\omega(s))_m - a_m f_m\right)^2}}{\|A\|_{min} - \|K\|_{op}} \leq \frac{\sqrt{\sum_{m=1}^\infty \left(\sum_{n=1}^\infty k_{mn} f_n + (\omega(s))_m - a_m f_m\right)^2}}{\alpha - \|k\|}.$$

This proves (4.21). □

In particular, if the kernel $k$ in Theorem 4.8 is a tensor product, we will have a corollary of Theorem 4.8.

**Corollary 4.9.** *Let $(S, \tau, \mu)$ be a topological probability space. Let $\{e_n\}_{n=1}^{\infty}$ be an orthonormal Schauder basis of $L_2[a,b]$. Let $g, h \in L_2[a,b]$ have sequences of coefficients $\{g_n\}_{n=1}^{\infty}$ and $\{h_n\}_{n=1}^{\infty}$, respectively.*

*Let $\omega: S \to L_2[a,b]$ be a single-valued mapping with sequence of coefficients $\{(\omega(s))_n\}_{n=1}^{\infty}$, for $s \in S$. Let $A = \{a_n\}_{n=1}^{\infty}$ be a bounded sequence of real numbers that induces a pointwise multiplication operator $A$. Let $\bar{s} \in S$. Suppose that $k, A$ and $\omega$ satisfy the following conditions.*

(a1)  $0 < \|g\|_2 \|h\|_2 < \|A\|_{\inf} \leq \|A\|_{op} \leq 1$;
(a2)  *The function* $s \to \omega(s)$ *is lower semicontinuous at* $\bar{s}$.

*Then, for any $\alpha$ with $\|A\|_{\inf} > \alpha > \|g\|_2\|h\|_2$, there exist a neighborhood $W_\alpha \subset S$ of $\bar{s}$ and a single-valued mapping $\sigma_\alpha: W_\alpha \to L_2[a,b]$ with sequence of coefficients $\{(\sigma_\alpha(s))_n\}_{n=1}^{\infty}$, for any $s \in S$ such that*

$$a_m(\sigma_\alpha(s))_m = g_m \sum_{n=1}^{\infty} h_n(\sigma_\alpha(s))_n + (\omega(s))_m, \text{ for any } s \in W_\alpha \text{ and for } m = 1, 2, \ldots. \quad (4.28)$$

*And, for any $f \in L_2[a,b]$, we have*

$$\sqrt{\sum_{m=1}^{\infty}((\sigma_\alpha(s))_m - f_m)^2} \leq \frac{\sqrt{\sum_{m=1}^{\infty}(\langle h,f\rangle g_m + (\omega(s))_m - a_m f_m)^2}}{\alpha - \|k\|}, \text{ for any } s \in W_\alpha.$$

*In particular, let $f = \theta$, that is, $f_m = 0$, for $m = 1, 2, \ldots$. We have*

$$\sqrt{\sum_{m=1}^{\infty}(\sigma_\alpha(s))_m^2} \leq \frac{\sqrt{\sum_{m=1}^{\infty}(\omega(s))_m^2}}{\alpha - \|k\|}, \text{ for any } s \in W_\alpha.$$

*Proof.* Let $k = g \times h$. By $\langle h, f \rangle = \sum_{n=1}^{\infty} h_n f_n$, we have

$$k(u,v) = g(u)h(v) = \left(\sum_{m=1}^{\infty} g_m e_m(u)\right)\left(\sum_{n=1}^{\infty} h_n e_n(v)\right)$$

$$= \sum_{m=1}^{\infty}\left(\sum_{n=1}^{\infty} g_m h_n e_m(u) e_n(v)\right), \text{ for every } (u,v) \in [a,b]^2.$$

This shows that, with respect to the orthonormal basis $\{e_m \times e_n\}_{m,n=1}^{\infty}$, the sequences of coefficients $\{k_{mn}\}_{m,n=1}^{\infty}$ of $k$ satisfies that $\{k_{mn}\}_{m,n=1}^{\infty} = \{g_m h_n\}_{m,n=1}^{\infty}$. In the representation of $k$ in Theorem 4.8, when $k_{mn}$ is substituted by $g_m h_n$, this corollary follows from Theorem 4.8 immediately. □

In particular in Corollary 4.9, if the considered bounded sequence $A$ has constant entries 1, then, the single-valued mapping $\sigma_\alpha: W_\alpha \to L_2[a,b]$ in Corollary 4.9 can be precisely solved.

**Corollary 4.10.** *Let $(S, \tau, \mu)$ be a topological probability space. Let $\{e_n\}_{n=1}^{\infty}$ be an orthonormal Schauder basis of $L_2[a,b]$. Let $g, h \in L_2[a,b]$ and $\omega: S \to L_2[a,b]$ be a single-valued mapping. Let $\bar{s} \in S$. Suppose that $k, A$ and $\omega$ satisfy the following conditions.*

(a1)  $0 < \|g\|_2 \|h\|_2 < 1$;
(a2)  *The function* $s \to \omega(s)$ *is lower semicontinuous at* $\bar{s}$.

*Then, for any $\alpha$ with $1 > \alpha > \|g\|_2\|h\|_2$, there exist a neighborhood $W_\alpha \subset S$ of $\bar{s}$ and a single-valued mapping $\sigma_\alpha: W_\alpha \to L_2[a,b]$ with sequence of coefficients $\{(\sigma_\alpha(s))_n\}_{n=1}^{\infty}$, for any $s \in S$ such that*

$$(\sigma_\alpha(s))_m = \frac{\langle \omega, h \rangle}{1 - \langle g, h \rangle} g_m + (\omega(s))_m, \text{ for any } s \in W_\alpha \text{ and for } m = 1, 2, \ldots. \quad (4.29)$$

*Proof.* Similarly, to the proof of Corollary 4.9, let $k = g \times h$ and let $A = \{1\}_{n=1}^{\infty}$ be the constant sequence

with all entries 1. Then, we have $0 < \|g\|_2\|h\|_2 < \|A\|_{\inf} = \|A\|_{op} = 1$. Hence, all conditions in Corollary 4.9 and in Corollary 4.3 are satisfied, which allows us to prove this corollary by the results of Corollary 4.9, or Corollary 4.3. Then, we have two ways to prove (4.29),

Way 1. By using (4.12) of Corollary 4.3 in the previous subsection. Let $\lambda = 1$ in (4.12), we obtain

$$(\sigma_\alpha(s))(u) = \frac{\int_a^b h(v)\omega(s)(v)dv}{1-\int_a^b g(v)h(v)dv} g(u) + \omega(s)(u), \text{ for any } u \in [a,b].$$

By (4.12), this is

$$(\sigma_\alpha(s))(u) = \frac{\langle \omega, h \rangle}{1-\langle g, h \rangle} g(u) + \omega(s)(u), \text{ for any } u \in [a,b].$$

This implies (4.29) immediately.

Way 2. Without using (4.12). We directly prove (4.29) by using (4.28). For the simplicity, we denote $(\sigma_\alpha(s))_m$ by $\sigma_m$ (That depends on $\alpha$ and $s$) and denote $(\omega(s))_m$ by $\omega_m$ (That depends on $s$). In (4.28), let $a_m = 1$, for all $m = 1, 2, \ldots$. Then, for any $s \in W_\alpha$, we have

$$\sigma_m = g_m \sum_{n=1}^\infty h_n \sigma_n + \omega_m, \text{ for } m = 1, 2, \ldots. \tag{4.30}$$

If $g_m = 0$, then $\sigma_m = \omega_m$. Hence, we assume $g_m \neq 0$. Write the equations (4.30) one by one and we obtain the following system of equations.

$$\left(h_1 - \frac{1}{g_1}\right)\sigma_1 + h_2\sigma_2 \quad + \quad h_3\sigma_3 \quad + \quad h_3\sigma_3 \quad + \cdots = -\frac{\omega_1}{g_1}$$

$$h_1\sigma_1 \quad + \left(h_2 - \frac{1}{g_2}\right)\sigma_2 + h_3\sigma_3 \quad + \quad h_4\sigma_4 \quad + \cdots = -\frac{\omega_2}{g_2}$$

$$h_1\sigma_1 \quad + \quad h_2\sigma_2 \quad + \left(h_3 - \frac{1}{g_3}\right)\sigma_3 + h_4\sigma_4 \quad + \cdots = -\frac{\omega_3}{g_3}$$

$$\ldots\ldots$$

By subtracting row $(m)$ − row $(1)$, for $m = 2, 3, \ldots$, we obtain

$$\left(h_1 - \frac{1}{g_1}\right)\sigma_1 + h_2\sigma_2 \quad + \quad h_3\sigma_3 \quad + \quad h_3\sigma_3 \quad + \cdots = -\frac{\omega_1}{g_1}$$

$$\frac{1}{g_1}\sigma_1 \quad \quad -\frac{1}{g_2}\sigma_2 \quad + \quad 0 \quad + \quad 0 \quad + \cdots = \frac{\omega_1}{g_1} - \frac{\omega_2}{g_2}$$

$$\frac{1}{g_1}\sigma_1 \quad + \quad 0 \quad \quad -\frac{1}{g_3}\sigma_3 \quad + \quad 0 \quad + \cdots = \frac{\omega_1}{g_1} - \frac{\omega_3}{g_3}$$

$$\ldots\ldots$$

Row $(1) + g_m h_m$ row $(m)$, for $m = 2, 3, \ldots$. The row (1) becomes

$$\left(h_1 - \frac{1}{g_1}\right)\sigma_1 + \frac{g_2 h_2}{g_1}\sigma_1 + \frac{g_3 h_3}{g_1}\sigma_1 + \cdots + 0 + \cdots = -\frac{\omega_1}{g_1} + g_2 h_2\left(\frac{\omega_1}{g_1} - \frac{\omega_2}{g_2}\right) + g_3 h_3\left(\frac{\omega_1}{g_1} - \frac{\omega_3}{g_3}\right) + \ldots$$

This infinitely dimensional system of linear equation becomes the following system

$$\frac{1}{g_1}(\langle g,h\rangle - 1)\sigma_1 = \frac{\omega_1}{g_1}(\langle g,h\rangle - 1) - \langle \omega, h\rangle.$$

$$\frac{1}{g_1}\sigma_1 \quad -\frac{1}{g_2}\sigma_2 \quad + 0 \quad + \quad 0 \quad + \cdots = \frac{\omega_1}{g_1} - \frac{\omega_2}{g_2}$$

$$\frac{1}{g_1}\sigma_1 \quad + \quad 0 \quad -\frac{1}{g_3}\sigma_3 \quad + \quad 0 \quad + \cdots = \frac{\omega_1}{g_1} - \frac{\omega_3}{g_3}$$

......

From the first equation, we have that $\sigma_1 = -\frac{\langle\omega,h\rangle}{\langle g,h\rangle - 1}g_1 + \omega_1$. Substituting this answer to row $(m)$, for $m = 2, 3, \ldots$, we have that $-\frac{1}{g_m}\sigma_m = \frac{\omega_1}{g_1} - \frac{\omega_m}{g_m} - \frac{1}{g_1}\left(-\frac{\langle\omega,h\rangle}{\langle g,h\rangle - 1}g_1 + \omega_1\right)$. This implies that

$$\sigma_m = -\frac{g_m\omega_1}{g_1} + \omega_m + \frac{g_m}{g_1}\left(-\frac{\langle\omega,h\rangle}{\langle g,h\rangle - 1}g_1 + \omega_1\right) = -\frac{\langle\omega,h\rangle}{\langle g,h\rangle - 1}g_m + \omega_m, \text{ for } m = 1, 2, 3, \ldots . \qquad \square$$

**Example 4.11.** Let $(S, \tau, \mu)$ be the standard topological probability space $[0, 1]$. Let $\{e_n\}_{n=1}^\infty$ be an orthonormal Schauder basis of $L_2[a, b]$. Let $g, h \in L_2[a, b]$ and $\omega: S \to L_2[a, b]$ respectively with their sequences of coefficients $\left\{\frac{1}{2^n}\right\}_{n=1}^\infty$, $\left\{\frac{1}{3^n}\right\}_{n=1}^\infty$ and $\left\{\frac{s^2}{4^n}\right\}_{n=1}^\infty$, for any $s \in S$. Let $A = \{1\}_{n=1}^\infty$ be the constant sequence with all entries 1. We see that $g, h, A$ and $\omega$ satisfy the following conditions.

(a1) $0 < \|g\|_2\|h\|_2 = \frac{1}{\sqrt{24}} < \|A\|_{\inf} = \|A\|_{\op} = 1$;

(a2) The function $s \to \omega(s)$ is continuous at every point $s \in S$.

Then, for any $\alpha$ with $1 > \alpha > \frac{1}{\sqrt{24}}$, there exists a single-valued mapping $\sigma_\alpha: [0, 1] \to L_2[a, b]$, that the sequence of coefficients $\{(\sigma_\alpha(s))_n\}_{n=1}^\infty$, for any $s \in S$ has the following explicit representation.

$$(\sigma_\alpha(s))_m = \left(\frac{5}{44} + \frac{1}{2^m}\right)\frac{s^2}{2^m}, \text{ for any } s \in [0, 1] \text{ and for } m = 1, 2, \ldots .$$

*Proof.* By Corollary 4.10, for any $s \in [0, 1]$ and for $m = 1, 2, \ldots$, we have

$$(\sigma_\alpha(s))_m = -\frac{\langle\omega,h\rangle}{\langle g,h\rangle - 1}g_m + (\omega(s))_m = -\frac{\frac{s^2}{11}}{\frac{1}{5}-1}\frac{1}{2^m} + \frac{s^2}{4^m} = \frac{\frac{s^2}{11}}{\frac{4}{5}}\frac{1}{2^m} + \frac{s^2}{4^m} = \left(\frac{5}{44} + \frac{1}{2^m}\right)\frac{s^2}{2^m}. \qquad \square$$

### 4.3. Parameterized Integral Equations

Let $k \in L_2[a, b]^2$. In contrast to the previous subsections, in this subsection, we particularly take the space $S = [a, b]$ and let $\omega: [a, b] \to L_2[a, b]$ be a measurable single-valued mapping (that is considered as a noise). Let $\bar{s} \in [a, b]$. For any given real number $\lambda$ with $|\lambda| \leq 1$, in this subsection, we will show that there is a neighborhood $W \subset [a, b]$ of $\bar{s}$ and a single-valued mapping $\sigma: W \to L_2[a, b]$ such that

$$\lambda(\sigma(s))(u) = \int_a^s k(u, v)(\sigma(s))(v)dv + \omega(s)(u), \text{ for any } s \in W \text{ and for any } u \in [a, b],$$

**Theorem 4.12.** *Let $k \in L_2[a, b]^2$ satisfy that there is $c \in (a, b)$ such that*

$$\left(\int_a^b \int_a^c |k(u, v)|^2 dv du\right)^{\frac{1}{2}} > 0. \tag{4.31}$$

Let $\omega: [a, b] \to L_2[a, b]$ be a single-valued mapping. Let $\lambda$ be a real number and $\bar{s} \in [c, b]$. Suppose that $k, \lambda$ and $\omega$ satisfy the following conditions.

(a1) $\|k\| < |\lambda| \leq 1$;
(a2) The function $s \to \omega(s)$ is lower semicontinuous at $\bar{s}$.

Then, for any $\alpha$ with $|\lambda| > \alpha > \|k\|$, there exist a neighborhood $W_\alpha \subset [c, b]$ of $\bar{s}$ and a single-valued mapping $\sigma_\alpha: W_\alpha \to L_2[a, b]$ such that

$$\lambda(\sigma_\alpha(s))(u) = \int_a^s k(u, v)(\sigma_\alpha(s))(v)dv + \omega(s)(u), \text{ for any } s \in W_\alpha \text{ and for any } u \in [a, b], \quad (4.32)$$

and $\quad \|\sigma_\alpha(s) - f\|_2 \leq \dfrac{\left(\int_a^b \left|\int_a^s k(u,v)f(v)dv + \omega(s)(u) - \lambda f(u)\right|^2 du\right)^{\frac{1}{2}}}{\alpha - \|k\|}, \text{ for any } s \in W_\alpha \text{ and } f \in L_2[a, b]. \quad (4.33)$

*Proof.* For any $s \in [c, b]$, define a single-valued function $\mathcal{X}_{[a,s]}$ on $[a, b]$ by

$$\mathcal{X}_{[a,s]}(v) = \begin{cases} 1, & \text{for } v \in [a, s], \\ 0, & \text{for } v \in (s, b], \end{cases} \text{ for any } v \in [a, b].$$

For any $s \in [c, b]$, define $k_s \in L_2[a, b]^2$ by $k_s(u, v) = k(u, v)\mathcal{X}_{[a,s]}(v)$, for any $(u, v) \in [a, b]^2$. Then, let $K_s$ be the integral operator corresponding to $k_s \in L_2[a, b]^2$. For any $f \in L_2[a, b]$, we have

$$K_s(f)(u) = \int_a^b k_s(u, v)f(v)dv = \int_a^b k(u, v)\mathcal{X}_{[a,s]}(v)f(v)dv = \int_a^s k(u, v)f(v)dv, \text{ for any } u \in [a, b]. \quad (4.34)$$

In particular, by the condition (4.31) in this theorem, at the point $c$, we have that $\|K_c\|_{op} = \left(\int_a^b \int_a^c |k(u, v)|^2 dv du\right)^{\frac{1}{2}} > 0$. By the conditions in this theorem, this implies that

$$0 < \|K_c\|_{op} \leq \|K_s\|_{op} \leq \|K\|_{op} \leq \|k\| < 1, \text{ for any } s \in [c, b]. \quad (4.35)$$

In AMZ Theorem, take $F = \lambda I$, in which $\lambda$ is the real number given in this theorem and $I$ is the identity mapping on $L_2[a, b]$. Meanwhile, let $P = [c, b]$ and define $G: L_2[a, b] \times [c, b] \to L_2[a, b]$ by

$$G(f, s)(u) = K_s(f)(u) + \omega(s)(u), \text{ for any } f \in L_2[a, b] \text{ and for any } s \in [c, b].$$

By (4.35), we can show that the mapping $G(\cdot, s): L_2[a, b] \to L_2[a, b]$ satisfies the Lipschitz condition on the whole space $L_2[a, b]$ relative to $L_2[a, b]$ for each $s \in [c, b]$ with the uniform modulus $\beta$ satisfying

$$0 < \|K_c\|_{op} \leq \beta \leq \|K\|_{op} \leq \|k\| < 1.$$

Then, all conditions (A1), (A2) and (A3) in Corollary 2.2 are satisfied, which proves this theorem. □

## 5. Conclusion

In the Arutyunov Mordukhovich Zhukovskiy Parameterized Coincidence Point Theorem, both of the involved mappings are set-valued mappings. In this paper, we consider two special cases of the Arutyunov Mordukhovich Zhukovskiy Theorem, which are presented as Corollaries 2.1 and 2.2 in Section 2. In Corollary 2.1, one of the considered mappings is single-valued mapping and the other one remains to be set-valued mapping. In this case, the existence of the coincidence point becomes an inclusion property. In Corollary 2.2, both of the involved mappings are single-valued mappings. Then, we get a parameterized coincidence point theorem for single-valued mappings, which is used in Section 4 to prove an existence theorem of solutions for some stochastic integral equations (See Theorem 4.1 and its corollaries).

Notice that the most difficult part for applying the Arutyunov Mordukhovich Zhukovskiy Parameterized Coincidence Point Theorem is the calculation of the covering constant for the involved set-valued mapping $F$. It is because to calculate the covering constant for $F$, one needs to calculate the Mordukhovich derivatives of $F$. From Section 2, we see that this is a very complicated procedure. This is why, in this paper, we only consider some linear and continuous mappings, for which the covering constant can be calculated. To extend the results in this paper to more general mappings, one needs to consider the following two questions.

**Question 1**. Do we find some practical and feasible techniques for calculating the Mordukhovich derivatives and the covering constants for set-valued mappings?

**Question 2**. Do we find some practical and feasible techniques for calculating the Mordukhovich derivatives and the covering constants for single-valued mappings without the Fréchet differentiability of the considered mappings?

**Author's declaration.** The author did not receive support from any organization for the submitted work. The author has no competing interests to declare that are relevant to the content of this article.
There is no any data supporting the results and analysis in the article.

**Acknowledgments** The author is grateful to Professor Robert Mendris, Professor Boris Mordukhovich, Professor Christiane Tammer and Professor Jen-Chih Yao for their communications and encouragements in the development stage of this paper. The author deeply thanks the anonymous referees for their valuable suggestions and helpful remarks, which allow me to improve the original manuscript.